\documentclass[nonblindrev]{informs3noheader}

\OneAndAHalfSpacedXI

\usepackage{mathtools}
\usepackage[inline]{enumitem}
\usepackage{subcaption}
\usepackage{booktabs}
\usepackage{array}
\usepackage{multirow}
\usepackage{booktabs,caption,fixltx2e}
\usepackage[flushleft]{threeparttable}
\usepackage{lscape}
\usepackage{makecell}
\usepackage{amssymb}
\usepackage{amsmath}
\usepackage[mathscr]{euscript}
\usepackage{algorithm}
\usepackage{algpseudocode}
\usepackage{pifont}
\usepackage{caption}
\usepackage{algorithmicx}
\usepackage{hyperref}
\usepackage{mathrsfs}
\usepackage{arydshln}
\usepackage{tabularx}
\usepackage{adjustbox}
\usepackage{bm, bbm}
\usepackage{soul}
\usepackage{blindtext}
\usepackage{longtable}
\hypersetup{
colorlinks=true,
linkcolor=black,
citecolor=black,
filecolor=black,
urlcolor=black,
}

\usepackage{comment}

\usepackage{tikz}
\usetikzlibrary{shapes.geometric,shapes.multipart}
\usetikzlibrary{positioning}
\usetikzlibrary{decorations.pathreplacing}
\definecolor{mygreen}    {RGB}{0,90,0}
\definecolor{myblue}     {RGB}{0,51,140}
\definecolor{myorange}   {RGB}{238,118,0}
\definecolor{myred}      {RGB}{126,0,0}
\definecolor{mygray}     {RGB}{100,100,105}
\definecolor{mygrayblue} {RGB}{0,128,128}
\definecolor{mygraygreen}{RGB}{128,128,0}
\definecolor{DarkPurple}     {RGB}{142, 36, 170}
\definecolor{LightPurple}    {RGB}{57, 130, 7}

\usepackage{natbib}
\bibpunct[, ]{(}{)}{,}{a}{}{,}%
\def\bibfont{\small}%

\usepackage[utf8]{inputenc}
\usepackage[english]{babel}
\usepackage{graphicx}



\newcolumntype{F}[1]{%
    >{\raggedright\arraybackslash\hspace{0pt}}p{#1}}%
\newcolumntype{T}[1]{%
    >{\centering\arraybackslash\hspace{0pt}}p{#1}}%

\TheoremsNumberedThrough
\EquationsNumberedThrough

\def\R{\mathbb{R}}

\def\calA{\mathcal{A}}
\def\calB{\mathcal{B}}
\def\calC{\mathcal{C}}

\def\calI{\mathcal{I}}
\def\calJ{\mathcal{J}}

\def\calN{\mathcal{N}}

\def\calP{\mathcal{P}}
\def\calQ{\mathcal{Q}}

\def\calT{\mathcal{T}}
\def\calU{\mathcal{U}}
\def\calV{\mathcal{V}}

\def\bi{\boldsymbol{i}}

\def\bx{\boldsymbol{x}}
\def\by{\boldsymbol{y}}
\def\bz{\boldsymbol{z}}

\def\bo{\boldsymbol{0}}

\def\st{\text{s.t.}}

\definecolor{mygreen}    {RGB}{0,90,0}
\definecolor{myblue}     {RGB}{0,51,140}
\definecolor{myorange}   {RGB}{238,118,0}
\definecolor{myred}      {RGB}{126,0,0}
\definecolor{mygray}     {RGB}{100,100,105}
\definecolor{mygrayblue} {RGB}{0,128,128}
\definecolor{mygraygreen}{RGB}{128,128,0}
\definecolor{DarkPurple}     {RGB}{142, 36, 170}
\definecolor{LightPurple}    {RGB}{57, 130, 7}

\begin{document}

\MANUSCRIPTNO{}

\ARTICLEAUTHORS{
	\AUTHOR{Adam Deng, Alexandre Jacquillat}
	\AFF{Massachusetts Institute of Technology,  Cambridge, MA}
}	

\RUNAUTHOR{Deng and Jacquillat}

\RUNTITLE{Coordinating mobile network coverage and vehicle routing}

\TITLE{Coordinating mobile network coverage and vehicle routing: a double column generation approach}

\ABSTRACT{The emergence of 5G technologies opens opportunities to support mission-critical activities with high-speed Internet coverage. This paper defines a joint job-emitting vehicle routing problem with time windows to coordinate the operations of mission-oriented vehicles (``mission vehicles'') and mobile emitters (``emitting vehicles''). This problem exhibits a joint vehicle routing structure, with coupling constraints to ensure that each job is supported by appropriate network coverage. We solve it via an exact and finite double column generation algorithm: pricing problems generate vehicle paths dynamically, and a master problem coordinates the operations of mission vehicles and emitting vehicles to ensure appropriate network coverage for each job. We propose several acceleration strategies to strengthen the algorithm's computational performance. Computational results show the scalability of the proposed methodology. Specifically, the methodology cuts runtimes by over 95\% in small-scale instances as compared to an explicit formulation, and scales to large-scale instances involving over 50 jobs. From a practical standpoint, results highlight the benefits of dynamically coordinating mission vehicles and emitting vehicles, thus suggesting opportunities to support emerging 5G technologies with dedicated routing algorithms.}


\maketitle

\vspace{-12pt}
\section{Introduction}

The advent of 5G technologies marks a major upgrade in the information and communication technology (ICT) infrastructure. Compared to its predecessors, the 5G network architecture provides three main service areas: (i) enhanced mobile broadband (eMBB) to power faster connectivity with lower latency; (ii) ultra reliable low latency communications (URLLC) to enable mission-critical applications requiring continuous and reliable coverage; and (iii) massive machine-type communications (mMTC) to power remote sensing capabilities by connecting multiple IoT devices in real time. These capabilities open new opportunities to support physical activities or perform them remotely, powered by high-speed connectivity. For example, police cars reacting to emergencies can upload files and gather data online before they reach the emergency scene; ambulances transporting patients in critical condition can spearhead treatment before reaching a healthcare facility by enabling remote interventions from a healthcare provider; military vehicles carrying soldiers and equipment can receive online updates and assessments through encrypted channels; other applications span smart agriculture, smart manufacturing, smart cities, etc.

\begin{figure}[h!]
    \centering
    \begin{subfigure}[t]{0.49\textwidth}
        \includegraphics[width=\textwidth]{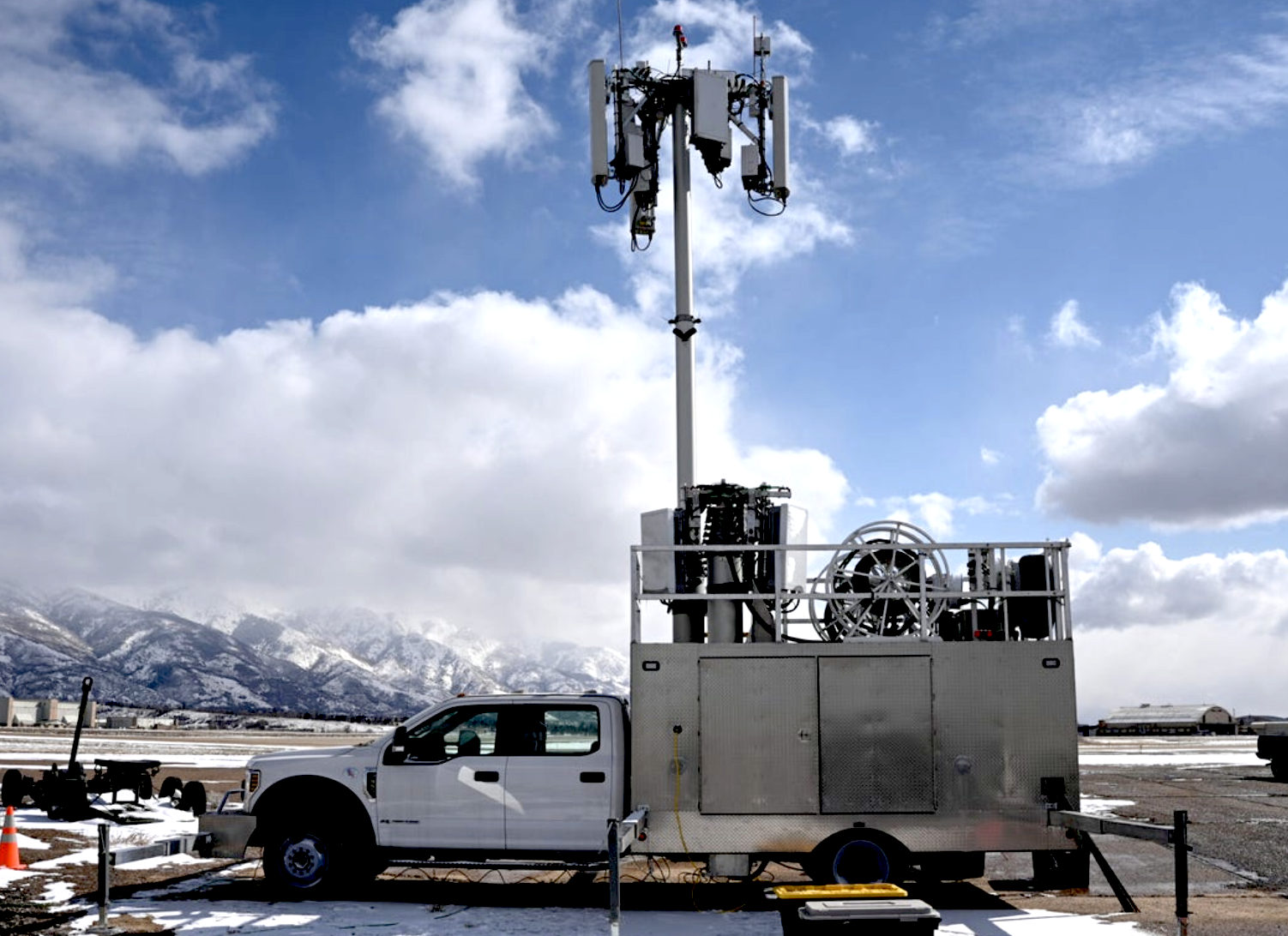}
        \caption{Truck with attached 5G mast}
        \label{fig:5gtruck}
    \end{subfigure}
    \begin{subfigure}[t]{0.49\textwidth}
        \includegraphics[width=\textwidth]{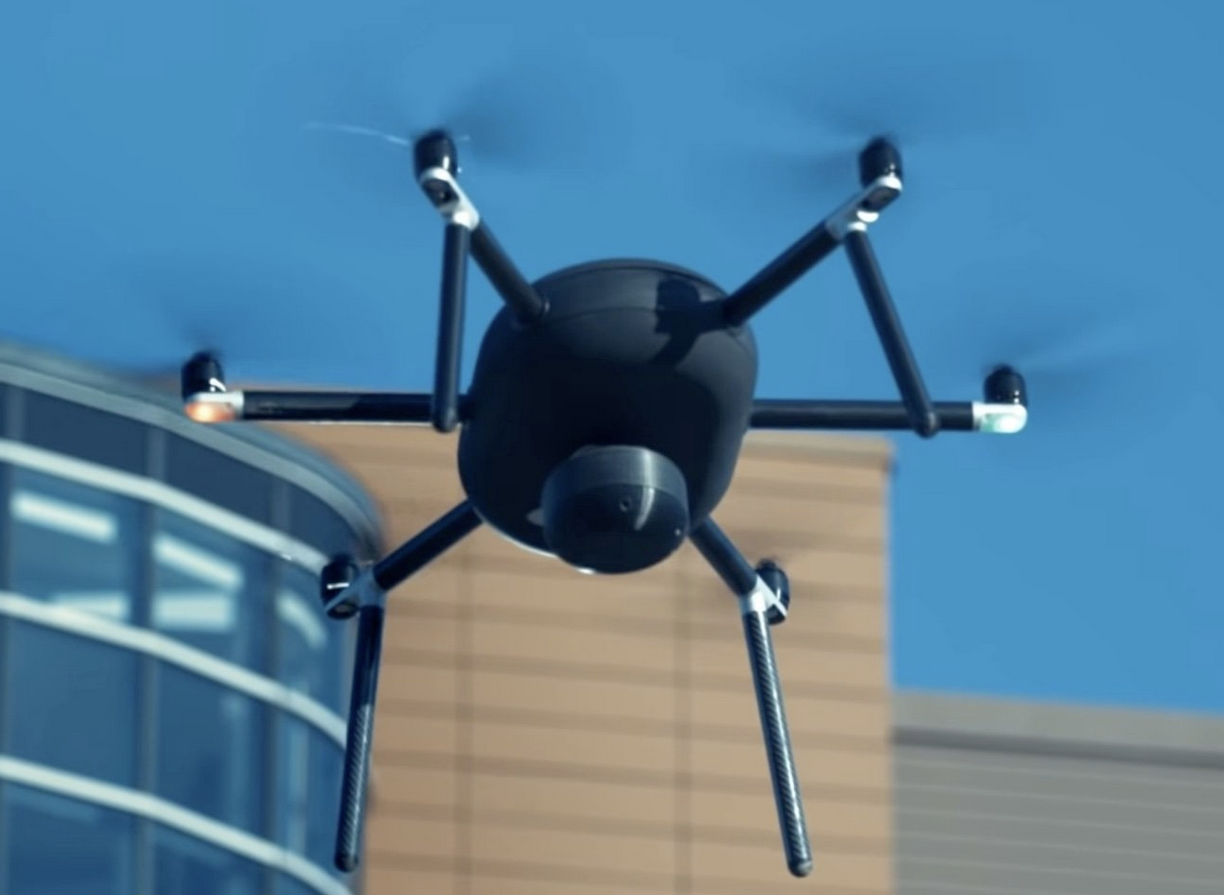}
        \caption{5G drone}
        \label{fig:movingrover}
    \end{subfigure}
    \caption{Prototypes of dynamic emitters as examples of the emitting vehicles considered in this paper.}
    \label{fig:emitters}
    \vspace{-12pt}
\end{figure}

These emerging mission-critical activities require high-speed local connectivity from the 5G infrastructure. One of the key elements of such ICT infrastructure is the set of mobile emitters that will be required to support missions spread across geographically disperse locations. Figure~\ref{fig:emitters} shows sample prototypes of mobile emitters, one involving a truck with an attached 5G mast or another one involving an emitter directly embedded in a drone. These options vary in terms of capital expenditures, operating expenditures, operational capabilities and coverage reliability, among other factors. Yet, all emitter prototypes evolve across a geographic area to provide local internet coverage. At the same time, these mobile emitters involve high acquisition costs and operating costs, along with significant operational complexity to ensure effective and reliable coverage for mission-critical activities. Currently, mobile emitters are not available for individuals, and 5G drones and trucks are still in the prototyping stage, including in Singapore. However, task-performing drones with 5G capabilities do exist, with the most advanced system the Qualcomm Flight RB5 5G Platform, costing over \$4,000, used for "mission-critical use cases" according to \href{https://www.modalai.com/pages/qualcomm-flight-rb5-5g-platform}{Modal AI}. To fully utilize the potential of high-speed coverage, emerging 5G technologies therefore require dedicated optimization and analytics capabilities to support mission-critical activities in the right place at the right time with a cost-effective mobile 5G infrastructure.

To address this question, this paper conceptualizes a new coordinated vehicle routing problem between mission-oriented vehicles (referred to as ``mission vehicles'') and mobile emitting devices (referred to as ``emitting vehicles''). This problem jointly optimizes the operations of emitting vehicles and mission vehicles. Both types of vehicles are coupled by the requirement that sufficient coverage must be available whenever a mission is being performed. Throughout this paper, we consider a closed system in an two-dimensional environment, in which the emitting vehicles are the sole sources of network coverage and mission vehicles are the sole users of the available 5G coverage. From a practical standpoint, this situation would correspond, for instance, to a military fleet performing a set of missions or search-and-rescue operations in a remote area, as opposed to a fleet of police vehicles or ambulances performing a set of jobs in a dense urban area. We leave these broader instances with multiple emitting sources and external users for future research.

The first contribution of this paper is to formulate an integrated vehicle routing problem that synchronizes the operations of mission vehicles and those of emitting vehicles (Section~\ref{sec:model}). Mission vehicles need to perform a set of jobs within corresponding time windows. Each job needs to be supported by local network coverage from the mobile emitters. The operations of the mobile emitters therefore determine which locations get connectivity at what times. A linking constraint ensures that each job gets performed during a continuous period of time while being consistently covered by the emitting vehicles. This problem features a coordinated vehicle routing structure with time windows, with a particular structure that emanates from the many-to-many mapping from emitting vehicles to mission vehicles---each emitting vehicle provides network coverage within its vicinity, and can therefore support multiple jobs simultaneously. The model's objectives include minimizing the number of emitting vehicles---motivated by their high acquisition costs---and minimizing the distance traveled by all vehicles---as a proxy toward minimizing operating costs.

We represent the problem in time-space networks to capture spatio-temporal coordination requirements between mission vehicles and emitting vehicles. We develop two formulations for the problem: (i) an explicit formulation with arc-based variables tracking operations in time-space networks, along with explicit constraints linking all vehicles' operations; and (ii) a set-partitioning formulation with path-based variables capturing the full operations of mission vehicles and emitting vehicles, along with a linking constraint ensuring that each job is performed when local network coverage is available. These two formulations are equivalent but the set-partitioning formulation defines a tighter continuous relaxation, albeit with an exponential number of path-based variables.

Our second contribution is to develop a double column generation algorithm to solve the set-partitioning formulation, with tailored acceleration strategies exploiting the coordination requirements between mission vehicles and coverage vehicles (Section~\ref{sec:algorithm}). The algorithm iterates between a master problem and independent pricing problems, and converges to an optimal solution in a finite number of iterations. Specifically, the master problem selects vehicle paths out of a restricted set of candidates, while enforcing the linking constraints to ensure that all mission activities are supported by local network coverage. The pricing problems generate new paths of negative reduced cost for mission vehicles and for coverage vehicles independently, or prove that none exists. The final step in the algorithm is to restore integrality, which we do through a stem and blender approach to augment the sets of mission and emission routes from the linear relaxation solution.

The success of the double column generation algorithm hinges on its ability to solve the pricing problems efficiently. Both pricing problems exhibit shortest path structures in time-space networks, which we solved using label-setting algorithms. The pricing problem for mission vehicles is naturally constrained by the time windows, which enable extensive tree pruning. The pricing problem for coverage vehicles, however, is much less constrained and can become much more time consuming. Our acceleration strategies restrict the shortest path tree for emitting vehicles by exploiting the structure of the coordinated routing problem. Specifically, we propose acceleration techniques based on (i) static time windows (emitting vehicles only need to visit locations when missions need to be performed in their vicinity); (ii) dynamic primal time windows (emitting vehicles need to visit locations where missions are being performed from the incumbent master problem solution); and (iii) dynamic dual time windows (an additional emission route is only considered if it loosens a binding coverage constraint). 

The third contribution of this paper is to demonstrate the scalability and practical benefits of the optimization methodology (Section~\ref{sec:results}). We develop an extensive experimental setup in collaboration with Singapore's Defence Science and Technology Agency (DSTA), using real-world data on the coverage of mobile emitters. From a technical standpoint, we show that the double column generation armed with our acceleration strategies scales to large and realistic instances of the problem. Specifically, we find that the explicit formulation barely scales to small and simplified instances (e.g., 10-30 missions and no coordination with emitter vehicles); baseline implementations of the double column generation algorithms scales to medium-sized instances of the problem (e.g., up to 35 jobs); in comparison, our acceleration strategies yield optimal solutions in larger instances, with over 50 jobs. We leverage the methodology to conduct a case study on the performance of coupled ICT-physical systems in support of mission-critical activities. Methodologically, our results show that double column generation with warm starts and all acceleration techniques is the optimal configuration for solving the most instances of the problem the quickest, and the number of jobs is the key variable in determining algorithm runtime, while the ratio of mission cluster radius versus coverage radius also plays a role (the smaller this ratio, the less time spent). Practically, we demonstrate that the optimized solution can provide significant benefits as compared to baseline solutions based on greedy heuristics and simple decomposition methods. We also underscore the benefits of mobile emitters by showing that the fleet of emitting vehicles can be significantly reduced if they are allowed to travel alongside mission vehicles instead of being stationary.

\section{Literature Review}

Since its introduction by \cite{dan59}, the vehicle routing problem (VRP) has become one of the most prominent combinatorial optimization problems in the operations research literature \citep[see][for detailed reviews]{fisher1995vehicle,golden2008vehicle,toth2014vehicle}. In its simplest form, the problem involves optimizing vehicle routes, each stating and ending in a depot, to serve a give set of customers to minimum total distance traveled. The canonical VRP has been extended to capture many practical considerations. For instance, the capacitated vehicle routing problem (VRP) considers a set of customers with fixed demand quantities, and imposes that each vehicle route can only serve demand up to a pre-determined capacity \citep{labbe1991capacitated,ralphs2003capacitated}. More closely related to our paper, the vehicle routing with time windows (VRPTW) optimizes vehicle routes subject to the additional constraint that each customer must be visited between acceptable times \citep{kolen1987vehicle}. Other extensions include vehicle routing with consistency requirements \cite{francis2008period,wang2022routing}, with split deliveries \citep{desaulniers2010branch}, with electric vehicles \cite{desaulniers2016exact}, with stochastic demand \citep{bertsimas1992vehicle}, with stochastic travel times \citep{laporte1992vehicle}, etc. In our paper, mission vehicles follow a VRPTW structure to perform jobs within time windows. However, the problem also departs from the literature in two ways. First, each job requires a certain amount of work, which creates a joint routing-scheduling structure in our problem. Second, and most importantly, each job needs to be performed while receiving network coverage from the emitting vehicles, thereby introducing coupling constraints between mission vehicles' routes and emitting vehicles' routes.

This distinction connects our paper to the fast-growing literature on coordinated routing, motivated by increasingly complex logistics operations in practice \citep[see][for a review]{Sav16}. One branch of this literature considers operations with flexible stopping locations. For instance, \cite{stieber2015multiple} and \cite{ozbaygin2017branch} optimize VRP operations when customers can be served in different places at different times of the day (e.g., home, work). \cite{gambella2018vehicle} and \cite{zhang2023routing} optimize VRP operations when customers can walk to the pickup locations or from the dropoff location, thereby also providing flexibility regarding the sequence of stops followed by the vehicles. The routes of the emitting vehicles in our paper exhibits similar flexibility. However, our problem features a different structure due to the impact of the stopping locations of emitting vehicles on the resulting network coverage for the mission vehicles. This problem creates complex coordinated routing requirements between two types of vehicles. A related problem is the horsefly routing problem, in which a primary vehicle (e.g., a ship, a truck) serves as a base for secondary vehicles (e.g., drones) \citep{carlsson2017coordinated,poikonen2019mothership}. Our setting departs from this problem because of the many-to-many mapping between mission vehicles and emitting vehicles (each emitting vehicle provides can support multiple mission vehicles and each mission vehicle can be supported by multiple emitting vehicles), as opposed to the one-to-many mapping between vehicles in earlier problems (e.g., a truck can support multiple drones but a drone needs to depart from one specific truck). As a consequence, the coordination between emitting vehicles and mission vehicles does not require them to be in the exact same location, because each emitting vehicle provides network coverage within its vicinity.

From a methodological standpoint, vehicle routing problems have been solved via branch-and-cut algorithms \citep{bard2002branch, ropke2007models}, set partitioning methods \citep{desaulniers2010branch,dabia2013branch} and time-space methods \citep{chardaire2005solving,baldacci2012new,dash2012time}. Due to the complexity of the problem, several metaheuristics have also been proposed, based notably on genetic algorithms \citep{cattaruzza2016multi, gansterer2018centralized}, simulated annealing \citep{Acc21}, tabu search \citep{xu2017truthful} and large-scale neighborhood search \citep{grangier2016adaptive}. In our paper, we adopt a column generation approach on time-space networks; however, to capture the coordination requirements between both types of vehicles, our methodology exhibits a double set-partitioning model armed with a double column-generation algorithm linking the routes of mission vehicles and of emitting vehicles. In particular, our algorithm features novel acceleration strategies to link the timing and coordination requirements of both types of vehicles in the (separate) pricing problems.

Finally, column generation is typically embedded into branch-and-price algorithms to restore integrality of the decision variables \citep{barnhart1998branch}. One of the challenges in branch-and-price, however, lies in the complex implementation requirements along with the extensive branching schemes that may be required to converge to an optimal integral solution. In this paper, we instead propose a stem-and-blender technique that creates perturbations along the optimal routes from the linear relaxation, and solves a restricted optimization problem on the resulting routes. Our computational results in this paper suggests that this approach can lead to optimal, or near-optimal solutions without incurring the computational costs of branch-and-price algorithms.

In summary, our paper contributes a new problem in coordinated vehicle routing that jointly optimizes the routes of emitting vehicles and of mission vehicles, subject to linking network coverage constraints. To solve it, it introduces a double column generation algorithm along with primal and dual acceleration strategies in the pricing problem as well as new integrality-forcing heuristics. As such, our paper contributes to the literature on vehicle routing and smart city operations to coordinate physical operations and ICT management in emerging mission-critical environments.

\section{Model formulation}
\label{sec:model}

\subsection{Problem statement}

Our problem involves dispatching mission vehicles to perform a set of jobs within pre-specified time windows as well as dispatching emitting vehicles to cover the mission vehicles while performing each job. We consider an area that constitutes a closed ICT system, meaning that the emitting vehicles are the sole sources of network coverage and the mission vehicles are the sole users. This setting captures coupled physical-ICT operations, such as search-and-rescue operations, wildfire suppression and military missions in remote areas that are not covered by high-speed broadband.

Let $U$ denote the number of mission vehicles available. In total, mission vehicles need to perform a set of $n$ jobs, stored in $\calI=\{1,\cdots,I\}$. Each job $i\in\calI$ needs to be performed for a duration $\tau_i$ within a time window characterized by the earliest start time $T_i^S$ and the latest end time $T_i^E$. Each mission vehicle starts and ends at the same depot; we use $i=0$ and $i=I+1$ to refer to the depot at the start and end of each mission vehicle's route, respectively. We assume that each mission vehicle needs to be connected to high-speed internet for the entire duration of the job. In contrast, we assume that mission vehicles do not need to be covered while traveling from one job to another, although this could be included in an extension of the proposed methodology. We also assume that each job needs to be performed continuously, meaning that a mission vehicle cannot start working on a job, stop because of a lack of internet connectivity, and resume working on the job thereafter.

Let $V$ denote the number of emitting vehicles available. Emitting vehicles are evolving in the same area as mission vehicles, and start and end in the same depot. We make no assumption on the speed of emitting vehicles; for instance, a truck with a 5G mast would travel at similar speeds as mission vehicles, whereas a 5G drone would likely travel faster. In principle, emitting vehicles could stop in any location in the map. For tractability purposes, we consider a set of $J$ locations in which emitting vehicles can stop, which we denote by $\calJ=\{1,\cdots,J\}$. Again, we use $j=0$ and $j=J+1$ to refer to the depot at the start and the end of each emitting vehicle's route. Each emitting vehicle provides network coverage, but the quality of coverage decreases as the distance between the emitter and the receiver increases. We model local network coverage via a binary many-to-many mapping, by assuming that each emitting vehicle provides network coverage to all mission vehicles located within a certain coverage radius. Specifically, we denote by $\calC_i\subseteq\calJ$ the subset of locations of emitting vehicles that cover each job location $i\in\calI$. Experimental data collected by prototype 5G emitters in Singapore suggest that high-speed connectivity is typically available within radii of 10--70 meters, with download speeds of 115--520 mbps (see Section~\ref{sec:results}). For reference, these speeds range from a fast home WiFi router and a Cat 6 ethernet cable at home.

\begin{figure}[h]
    \centering
    \includegraphics[width=0.6\textwidth]{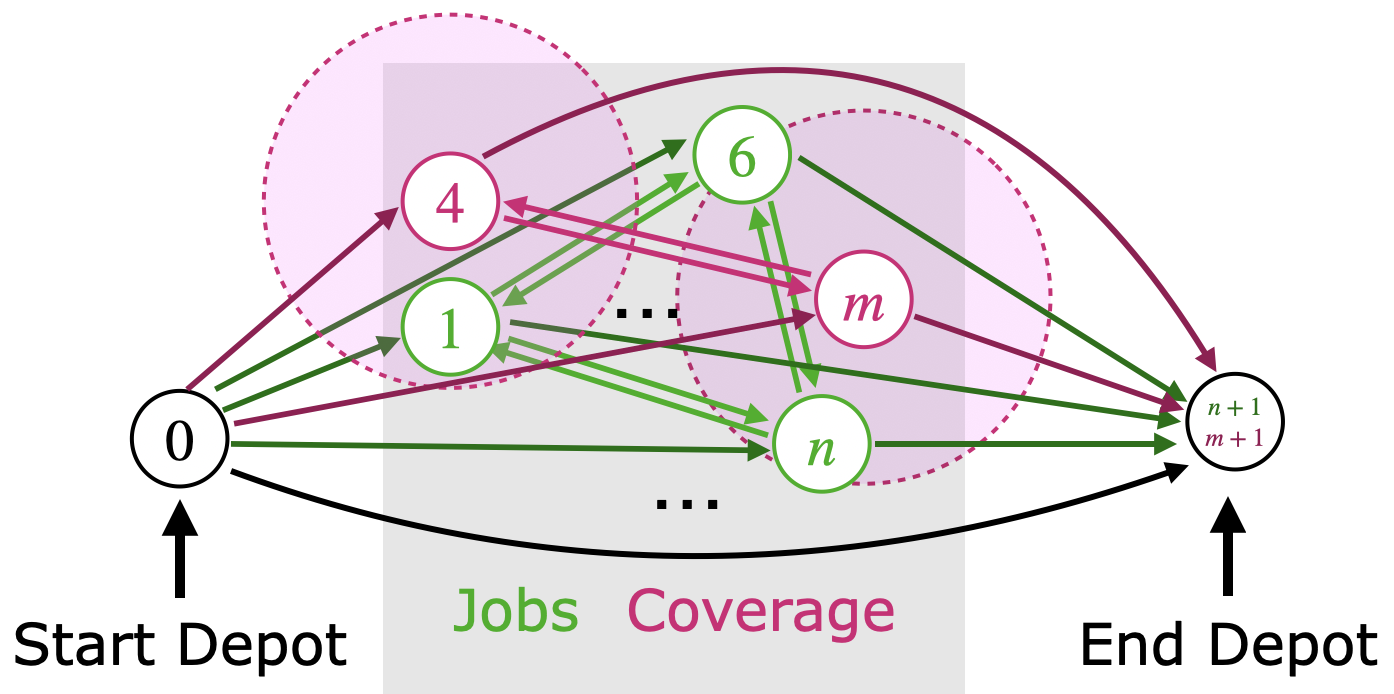}
    \caption{Graph of mission vehicles' operations (green nodes and arcs) and emitting vehicles' operations (purple nodes and arcs). Emitting vehicles provide local network coverage within a radius (purple circles).}
    \label{fig:graph}
\end{figure}

Figure~\ref{fig:graph} represents the physical operations of the problem. Note that the structure gives rise to spatial-temporal coordination requirements between mission vehicles and emitting vehicles, due to the fact that jobs can only performed in locations and at times that are covered by the emitting vehicles. Within the literature on coordinated vehicle routing, the problem distinguishes itself due to the many-to-many mapping structure between mission vehicles and emitting vehicles---each mission vehicle can be supported by any emitting vehicle, and any emitting vehicle can support all mission vehicles located within the coverage radius.

In order to capture the coordination requirements between mission vehicles and emitting vehicles, we represent operations in time-space networks. This representation relies on a time discretization into a set $\calT$, where $t=0$ and $t=T+1$ corresponding to the start and end of the planning horizon, respectively. Specifically, we create one time-space network for mission vehicles and one time-space network for emitting vehicles; this representation relies on the assumption that all mission vehicles are homogeneous and all emitting vehicles are homogeneous, but this assumption can be easily relaxed in an extension of the model. In this representation, each node represents a location-time pair, and each arc captures a vehicle trip from one location to another one. This representation integrates vehicle speeds into the design of time-space arcs. Note that the set of nodes and the set of arcs can vary across time-space networks, reflecting the different set of stopping locations for mission vehicles and for emitting vehicles as well as potentially different speeds. By design, flow balance in these time-space networks will capture continuity of operations in time and space.

Our formulation minimizes the distance traveled by mission vehicles and emitting vehicles. We encapsulate the distance via arc-based parameters $\{c^M_a:a\in\calA^M\}$ and $\{c^E_a:a\in\calA^E\}$ for mission vehicles and emitting vehicles, respectively; for future reference, we also denote by $\{d^M_{ij}:i,j\in\calI\}$ and $\{d^E_{ij}:i,j\in\calJ\}$ the corresponding distances in the physical network. This distance objective serves as a proxy for operating costs. Note that, by varying the number of available mission vehicles $B$ and the number of available emitting vehicles $V$, the model formulates a bi-objective problem, by trading off operating costs (distance traveled) and acquisition costs (fleet size). In particular, the model can also be leveraged to determine the smallest possible fleet of emitting vehicles to perform a set of missions. We report results based on both objectives in Section~\ref{sec:results}.

\subsection{Explicit optimization formulation}

\begin{table}[h]
    \centering
\caption{Notation for explicit formulation.}\label{tab:notation}
{\footnotesize \begin{tabular}{llp{13.5cm}}
\toprule[1pt]
    \bf Notation & \bf Type & \bf Description \\ \hline
    $\calI$ & Set & Job, or job locations \\
    $\calJ$  & Set & Coverage locations for emitting vehicles\\
    $\calC_i$  & Set & Subset of locations of emitting vehicles that cover job $i\in\calI$\\
    $\calT$ & Set & Time periods\\
    $\calU$ & Set & Set of mission vehicles\\
    $\calV$ & Set & Set of emitting vehicles\\
    $\calN^M$ & Set & Time-space nodes for mission vehicles: $n=(i,t)\in\calN^M$ if $i\in\calI$ and $t\in\calT$\\
    $\calN^E$ & Set & Time-space nodes for emitting vehicles: $n=(j,t)\in\calN^E$ if $j\in\calJ$ and $t\in\calT$\\
    $\calA^M$ & Set & Time-space arcs for mission vehicles: $(n,m)\in\calA^M$ if a vehicle can travel from $n\in\calN^M$ to $m\in\calN^M$\\
    $\calA_n^{M,-}$ & Set & Time-space arcs for mission vehicles that are incoming into node $n\in\calN^M$\\
    $\calA_n^{M,+}$ & Set & Time-space arcs for mission vehicles that are outgoing from node $n\in\calN^M$\\
    $\widetilde\calA_i^M$ & Set & Time-space arcs for mission vehicles that are idle in location $i\in\calI$\\
    &&$(n,m)\in\widetilde\calA_i^M\iff n,m\in\calN^M,\ n=(i,t),\ m=(i,t+1)$ for some $t\in\calT$\\
    $\calA^E$ & Set & Time-space arcs for emitting vehicles: $(n,m)\in\calA^E$ if a vehicle can travel from $n\in\calN^E$ to $m\in\calN^E$\\
    $\calA_n^{E,-}$ & Set & Time-space arcs for emitting vehicles that are incoming into node $n\in\calN^E$\\
    $\calA_n^{E,+}$ & Set & Time-space arcs for emitting vehicles that are outgoing from node $n\in\calN^E$\\
    $\widetilde\calA_j^E$ & Set & Time-space arcs for emitting vehicles that are idle in location $j\in\calJ$\\
    &&$(n,m)\in\widetilde\calA_i^E\iff n,m\in\calN^E,\ n=(j,t),\ m=(j,t+1)$ for some $t\in\calT$\\
    $c^M_a$ & Parameter & Cost of traveling along arc $a\in\calA^M$ for a mission vehicle\\
    $c^E_a$ & Parameter & Cost of traveling along arc $a\in\calA^E$ for an emitting vehicle\\
    $T_i^S$ & Parameter & Start of time window for job $i\in\calI$ \\
    $T_i^E$ & Parameter & End of time window for job $i\in\calI$ \\
    $\tau_i$ & Parameter & Work load of job $i\in\calI$ \\
    \bottomrule[1pt]
\end{tabular}}
\end{table}

Table~\ref{tab:notation} summarizes the model's notation. We define the following decision variables.
\begin{align*}
    z_{ua} = & \begin{cases}1&\text{if mission vehicle $u\in\calU$ travels along time-space arc $a\in\calA^M$}\\0&\text{otherwise}\end{cases} \\
    x_{va} = & \begin{cases}1&\text{if emitting vehicle $v\in\calV$ travels along time-space arc $a\in\calA^E$}\\0&\text{otherwise}\end{cases} \\
    y_{it}^S= &  \begin{cases}1&\text{if job $i\in\calI$ has started by time $t\in\calT$}\\0&\text{otherwise}\end{cases} \\
    y_{it}^E= &  \begin{cases}1&\text{if job $i\in\calI$ has ended by time $t\in\calT$}\\0&\text{otherwise}\end{cases}
\end{align*}

The variables $z_{ua}$ and $x_{va}$ reflect the routing of the mobile vehicles and emitting vehicles in their corresponding time-space networks, respectively. The variables $y_{it}^S$ and $y_{it}^E$ capture the start time and the end time of each job. Note that each row of these matrices encodes the time \textit{by} which a job starts or ends, as opposed to the time \textit{at} which a job starts or ends. In other words, each row includes a sequence of zeros followed by a sequence of ones, where the first ``1'' denotes the time at which a job starts or ends. For instance, if $y^S_{1,\cdot}=(0,0,1,1,1,\cdots)$ and $y^E_{1,\cdot}=(0,0,0,0,1,\cdots)$, then job $1$ starts in period 3 and ends in period 5. The structure of these variables is inspired by \cite{bertsimas1998air} in air traffic management, and strengthens the model's structure.

The optimization formulation is given as follows:
\begin{align}
\min  &\quad \displaystyle  \sum_{u \in \calU} \sum_{a \in \calA^M} c^M_az_{ua} + \sum_{v \in \calV} \sum_{a \in \calA^E} c^E_ax_{va} \label{df_objective}\\
    	\text{s.t} 	&\quad\textbf{Constraints for mission vehicles}\nonumber\\
    	&\quad \displaystyle \sum_{t \in \calT} \sum_{a \in\calA^{M,+}_{(0, t)}} z_{ua} = 1, \ \forall u \in \calU \label{jobleaveonce} \\
    	&\quad\displaystyle  \sum_{t \in \calT} \sum_{a \in \calA^{M,-}_{(I+1, t)}} z_{ua} = 1, \ \forall u \in \calU \label{jobenteronce}\\
    	&\quad\displaystyle \sum_{a \in \calA^{M,-}_n} z_{ua} = \sum_{a \in \calA^{M,+}_n} z_{ua}, \ \forall u \in \calU,\ \forall n \in \calN^M \label{conserveflowjob}\\
         &\quad\displaystyle \sum_{u \in \calU} \sum_{t\in\calT}\sum_{a \in \calA^{M,-}_{(i,t)}} z_{ua} = 1, \ \forall i \in \calI
         \label{jobpartition}\\
    	&\quad\textbf{Constraints for emitting vehicles}\nonumber\\
    	&\quad \displaystyle \sum_{t \in \calT} \sum_{a \in\calA^{E,+}_{(0, t)}} x_{va} = 1, \ \forall v \in \calV \label{coverageleaveonce} \\
    	&\quad\displaystyle  \sum_{t \in \calT} \sum_{a \in \calA^{E,-}_{(J+1, t)}} x_{va} = 1, \ \forall v \in \calV \label{coverageenteronce}\\
    	&\quad\displaystyle \sum_{a \in \calA^{E,-}_n} x_{va} = \sum_{a \in \calA^{E,+}_n} x_{va} \ \forall u \in \calU,\ \forall n \in \calN^E, \label{conserveflowcpv}\\
    	&\quad\textbf{Time window constraints}\nonumber\\
            &\quad\displaystyle y_{it}^S \leq y_{i, t+1}^S, \ \forall i \in \mathcal{I},\ \forall t \in \calT \label{timestart} \\
            &\quad\displaystyle y_{it}^E \leq y_{i, t+1}^E, \ \forall i \in \mathcal{I},\ \forall t \in \calT \label{timeend} \\
            &\quad\displaystyle y_{it}^S = 0, \ \forall t < T_i^S, \ \forall i \in\calI \label{notstarted} \\
            &\quad\displaystyle y_{it}^E = 1, \ \forall t \geq T_i^E, \ \forall i \in\calI \label{ended} \\
            &\quad\displaystyle \sum_{t \in \calT} (y_{it}^S - y_{it}^E) \geq \tau_i, \ \forall i \in\calI \label{finishjob} \\
            &\quad\textbf{Linking constraints}\nonumber\\
        	&\quad\displaystyle y_{it}^S - y_{it}^E \leq \sum_{u \in \calU}\ \sum_{a\in\widetilde\calA_i^M} z_{ua}, \ \forall i \in\calI,\ \forall t \in \calT\label{workcannotbedone}\\
            &\quad\displaystyle y_{it}^S - y_{it}^E \leq \sum_{v \in \calV}\ \sum_{a\in\widetilde\calA_j^E:j\in \calC_i} x_{va}, \ \forall i\in\calI,\ \forall t\in\calT\label{sufficientcov}\\
        	&\quad\textbf{Domain of definition:}\ \text{$\bx,\ \bz,\ \by^S,\ \by^E$ binary} \label{domain}
\end{align} 

Equation~\eqref{df_objective} minimizes the total operating cost across mission vehicles and emitting vehicles. Constraints~\eqref{jobleaveonce} and~\eqref{jobenteronce} guarantee that each mission vehicle leaves the depot and returns to the depot exactly once. Constraint~\eqref{conserveflowjob} ensures flow conservation for mission vehicles, by stating that any incoming arc must be accompanied by an outgoing arc at every node. Equations~\eqref{coverageleaveonce}--\eqref{conserveflowcpv} enforce analogous flow constraints for emitting vehicles. Constraint~\eqref{jobpartition} guarantees that each job is performed by exactly one mission vehicle. Equations~\eqref{timestart} and~\eqref{timeend} enforce the monotonicity of the start time and end time variables, consistently with their definition. Constraints~\eqref{notstarted} and~\eqref{ended} impose the time windows. Constraint~\eqref{finishjob} enforces the workload required to complete each job.

The remaining constraints are linking constraints between all variables. The first one ensures consistency between the operations of mobile vehicles and the completion of the jobs: Constraint~\eqref{workcannotbedone} ensures that work can only be performed at a job if at least one mobile vehicle is idle at the corresponding location. The second one ensures consistency between the operations of mobile vehicles and of emitting vehicles: Constraint~\eqref{sufficientcov} states that work can only be performed in job $i$ at time $t$ (i.e., $y_{it}^S - y_{it}^E=1$) if an emitting vehicle is present in at least one of its covering locations.

Whereas the explicit model provides a natural framework for the coordinated routing problem, it defines a weak optimization formulation. First, the number of variables, albeit polynomial, grows quadratically with the spatial-temporal discretization of the network. Second, the formulation defines a weak linear optimization relaxation due to the linking constraints between workload variables and the operations of the mobile vehicles---which perform the jobs---as well as those of emitting vehicles---which define network coverage. As such, the problem complicates an already-complex vehicle routing structure with time windows in two respects: (i) by separating the workload variables from the routing variables, and (ii) by applying the coordination requirements via linking constraints. As we shall see in our computational results, these complexities severely hinder the scalability of the explicit formulation, which can only handle small-sized and simplified instances.

\subsection{Set-partitioning formulation}

To address the difficulties of the explicit formulation, we propose a set-partitioning formulation based on path-based variables---for both types of vehicles. In addition to the notation defined in Table~\ref{tab:notation}, we define a set of paths $\calQ$ for mission vehicles and a set of paths $\calP$ for emitting vehicles. We introduce additional parameters to capture the cost of each path as the sum of the corresponding arc costs, and to map each path with each location-time pair. As such, the paths of the mission vehicles define vehicle availability to perform work in each job location, and the paths of the emitting vehicles define network coverage. Specifically, a path $q \in \calQ$ for mission vehicle $u\in\calU$ is defined by variables $z_{ua}$ for all $a\in\calA^M$ and $y_{it}$ for all $i\in\calI,t\in\calT$, which collectively satisfy Constraints~\eqref{jobleaveonce}--\eqref{jobpartition} and Constraints~\eqref{timestart}--\eqref{workcannotbedone}. Similarly, a path $p \in \calP$ for vehicle $v\in\calV$ is defined by variables $x_{va}$ for all $a\in\calA^E$ satisfying Constraints~\eqref{coverageleaveonce}--\eqref{conserveflowcpv}. Note, importantly, that the paths of mission vehicles are defined to encapsulate both vehicle routes in the physical network and work activities at the job locations. In other words, the set $\calQ$ contains several paths that visit the same set of job locations in the same sequence but perform work in one or multiple jobs at different times. The new notation for the set-partitioning formulation is provided in Table~\ref{tab:SPnotation}.

\begin{table}[h!]
    \centering
\caption{Additional notation for the set-partitioning formulation.}\label{tab:SPnotation}
{\footnotesize \begin{tabular}{llp{12cm}}
\toprule[1pt]
    \bf Notation & \bf Type & \bf Description \\ \hline
    $\calQ$ & Set & Paths for mission vehicles\\
    $\calP$  & Set & Paths for emitting vehicles\\
    $\calQ_i$ & Set & Subset of paths for mission vehicles that visit job $i\in\calI$\\
    $\overline{c}^M_q$ & Parameter & Cost of path $q\in\calQ$ for a mission vehicle\\
    $\overline{c}^E_p$ & Parameter & Cost of path $p\in\calP$ for an emitting vehicle\\
    $\delta_{it}^q$ & Parameter & Binary indicator of whether path $q\in\calQ$ performs work in job $i\in\calI$ at time $t\in\calT$ \\
    $g_{it}^p$ & Parameter & Binary indicator of whether path $p\in\calP$ covers job $i\in\calI$ at time $t\in\calT$\\
    \bottomrule[1pt]
\end{tabular}}
\end{table}

We introduce the following path-based decision variables. Note that it is sufficient to consider binary variables because each job is visited once by one vehicle, and because there is no value in multiple emitting vehicles traveling along the same path. Moreover, these variables do not need to be indexed by the vehicle, which further simplifies the model's structure by eliminating symmetry.
\begin{align*}
    \overline{z}_{q} = & \begin{cases}1&\text{if a mission vehicle travels along path $q\in\calQ$}\\0&\text{otherwise}\end{cases} \\
    \overline{x}_{p} = & \begin{cases}1&\text{if an emitting vehicle travels along path $p\in\calQ$}\\0&\text{otherwise}\end{cases}
\end{align*}

The set-partitioning formulation is given as follows. It still minimizes operating costs (Equation~\eqref{st_objective}), subject to set-partitioning constraints ensuring that each job is visited by exactly one path (Equation~\eqref{st_jobvisitedonce}), that the chosen paths do not exceed the fleet capacity (Equations~\eqref{st_atmostjobpaths} and~\eqref{st_atmostcovpaths}), and that mission vehicles are covered by emitting vehicles whenever performing work (Equation~\eqref{st_linker}). This model combines two path-based set-partitioning structures for mission vehicles and emitting vehicles, along with a linking constraint in Equation~\eqref{st_linker}.
\begin{align} 
\min  &\quad \displaystyle  \sum_{q \in \calQ} \overline{c}^M_q\overline{z}_{q}+\sum_{p \in \calP} \overline{c}^E_p\overline{x}_p \label{st_objective}\\
    	\text{s.t} 	
            &\quad \displaystyle \sum_{q \in \calQ_i} \overline{z}_{q} = 1, \ \forall i \in\calI \label{st_jobvisitedonce} \\
            &\quad\displaystyle \sum_{q \in \calQ} \overline{z}_{q} \leq U \label{st_atmostjobpaths}\\
            &\quad\displaystyle \sum_{p \in \calP} \overline{x}_{p} \leq V \label{st_atmostcovpaths}\\
    	&\quad\displaystyle \sum_{q \in \calQ} \delta_{it}^q\overline{z}_q \leq \sum_{p \in \calP} g_{it}^p\overline{x}_p, \ \forall i \in\calI,\ \forall t \in \calT \label{st_linker}\\
        &\quad\textbf{Domain of definition:}\ \text{$\overline\bx,\ \overline\bz$ integer} \label{SPdomain}
\end{align}

The set-partitioning formulation is equivalent to the explicit formulation. However, it eliminates several sets of constraints by embedding them into the definition of a path, including the flow-balance constraints (Equations~\eqref{jobleaveonce}--\eqref{conserveflowjob} and Equations~\eqref{coverageleaveonce}--\eqref{conserveflowcpv}), the time-window constraints (Equations~\eqref{timestart}--\eqref{finishjob}), and the linking constraint between the routing variables for mission variables and the workload variables. Obviously, the set-partitioning formulation comes at the cost of an exponential number of path-based variables. Recall, also that the proposed set-partitioning formulation relied on a definition of a path for mission vehicles that encapsulates both their route across job locations and their workload on the jobs. This structure comes at the cost of a larger number of path-based variables (each physical route is duplicated into several paths corresponding to different work plans) but strengthens the formulation by convexifying the time-window constraints and the routing-workload linking constraint. Proposition~\ref{prop:formulations} summarizes and proves these observations.

\begin{proposition}\label{prop:formulations}
    The set-partitioning formulation and the explicit formulation are equivalent, but the set-partitioning formulation defines a tighter linear relaxation.
\end{proposition}

To circumvent the exponential number of path-based variables in the set-partitioning formulation, we proceed via column generation in Section~\ref{sec:algorithm}. In fact, we propose a double column generation algorithm that separates the generation of paths for mission vehicles and for emitting vehicles, and brings them together in the master problem via the linking constraint. Before proceeding, we identify some structural properties of the optimal solution that will enable us to simplify the set of paths $\calQ$ for mission vehicles, and therefore strengthen the overall model and algorithm.

Proposition~\ref{prop:workwhenstay} states that we can restrict the paths of mission vehicles to those that perform work whenever idle in a job location. In other words, any slack time can be absorbed while traveling between locations. This property enables to circumvent the workload planning component of the problem to focus on the routing structure of the problem when generating paths for mission vehicles. This result reduces the size of the set-partitioning formulation and, most importantly, will be instrumental to accelerate the pricing problem in the column generation algorithm.

\begin{proposition}\label{prop:workwhenstay}
    In set-partitioning and column generation, it is without loss of generality to assume that mission vehicles only stay at job locations when performing work.
\end{proposition}

\section{A double column generation algorithm}
\label{sec:algorithm}

We propose a double column generation (DCG) algorithm to circumvent the exponential number of variables in the set-partitioning formulation. The formulation exhibits a vehicle routing structure with time windows for mission vehicles, and a vehicle routing structure with flexible locations for emitting vehicles. Together, they give rise to a coupled vehicle routing structure. Accordingly, we propose a decomposition scheme with a restricted master problem that selects paths for all vehicles while enforcing the linking constraints, and pricing problems that separate the operations of mission vehicles and emitting vehicles into independent shortest path formulations.

\subsection{Overview of algorithm}

The DCG algorithm maintains a subset $\calQ^0$ of paths for mission vehicles (similarly, $\calQ^0_i$ denotes the subset of paths that visit job $i\in\calI$) and a subset $\calP^0$ of paths for emitting vehicles. The Restricted Master Problem (RMP) applies the linear relaxation of the set-partitioning formulation based on these restricted path options. It is given as follows:
\begin{align} 
\min  &\quad \displaystyle  \sum_{q \in \calQ^0} \overline{c}^M_q\overline{z}_{q}+\sum_{p \in \calP^0} \overline{c}^E_p\overline{x}_p \label{dcg_obj}\\
    	\text{s.t} 	
            &\quad \displaystyle \sum_{q \in \calQ^0_i} \overline{z}_{q} = 1, \ \forall i \in\calI \label{dcg_jobvisitedonce} \\
            &\quad\displaystyle \sum_{q \in \calQ^0} \overline{z}_{q} \leq U \label{dcg_atmostjobroutes}\\
            &\quad\displaystyle \sum_{p \in \calP^0} \overline{x}_{p} \leq V \label{dcg_atmostcovroutes}\\
    	&\quad\displaystyle \sum_{q \in \calQ^0} \delta_{it}^q\overline{z}_q \leq \sum_{p \in \calP} g_{it}^p\overline{x}_p, \ \forall i \in\calI,\ \forall t \in \calT \label{dcg_linker}\\
        &\quad\textbf{Domain of definition:}\ \text{$\overline\bx\geq\bo,\ \overline\bz\geq\bo$} \label{dcgdomain}
\end{align}

Let $\pi_i\in\R$, $\rho\in\R_+$, $\beta\in\R_+$ and $\xi_{it}\in\R_+$ denote the dual variables associated with Equations~\eqref{dcg_obj}--\eqref{dcg_linker}, respectively. Per linear optimization theory, the RMP solution is optimal in the linear set-partitioning relaxation if the reduced costs of all (omitted) variables is non-negative, i.e.:
\begin{align}
    &\quad\displaystyle C^q + \rho - \sum_{i \in\calI:q\in\calQ_i} \pi_i + \sum_{i \in\calI} \sum_{t \in \calT} \delta_{it}^q \xi_{it}\geq0,\ \forall q \in \calQ \label{dcg_jobssubproblem} \\
    &\quad \displaystyle C^p + \beta - \sum_{i \in\calI} \sum_{t \in \calT} g_{it}^p \xi_{it}\geq0,\ \forall p \in \calP \label{dcg_covssubproblem}
\end{align}

Equations~\eqref{dcg_jobssubproblem} and~\eqref{dcg_covssubproblem} give rise to the pricing problem formulations in our DCG algorithm by seeking paths of negative reduced cost for mission vehicles and emitting vehicles, respectively, or prove that none exists. We develop algorithms to solve the pricing problem for mission vehicles in Section~\ref{subsec:mission}, and for emitting vehicles in Section~\ref{subsec:emitting}. We present the full algorithm in Section~\ref{sec:fullAlgorithm}.

\subsection{Pricing problem: seeking paths for mission vehicles}
\label{subsec:mission}

The pricing problem for mission vehicles seeks the path of minimal reduced cost, as follows:
\begin{align}
    \min &\quad\displaystyle \sum_{a\in\calA^M}c^M_az_{ua} + \rho - \sum_{i\in\calI}\sum_{t\in\calT}\sum_{a\in\calA_{(i,t)}^{M,-}} \pi_iz_{ua} + \sum_{i\in\calI}\sum_{t\in\calT} \xi_{it}y_{it}
    \label{jo_jobssubproblem}\\
    \st&\quad\text{Constraints~\eqref{jobleaveonce}--\eqref{jobpartition} and Constraints~\eqref{timestart}--\eqref{workcannotbedone}; $\bz,\by$ binary}
\end{align}

This problem features a prize-collecting shortest path structure with time windows. Each job $i\in\calI$ is associated with a reward $\pi_i$ included in any incoming arc. Each arc comes at a traveling cost $c_a^M$, which decomposes the total cost of the path $C^q$. Next, each vehicle that stays in job location $i\in\calI$ at time $t\in\calT$ incurs a cost $\xi_{it}$, reflecting the shadow cost of bringing an emitting vehicle in the vicinity of the job location to provide local network coverage. Finally, the fixed cost $\rho$ is added at the beginning of any path to capture the shadow cost of using one mission vehicle.

Recall from Proposition~\ref{prop:workwhenstay} that any mission vehicle stays at a job location only while performing work. This result enables us to eliminate the $\by$ variable from the formulation by adding constraints ensuring that the vehicle arrives at the job location after the start of the time window, that it leaves before the end of the time window, and that it stays for a duration exactly equal to the workload requirements. This property is critical to decompose the pricing problem into a shortest path path in the graph of job locations. Specifically, the problem can be reformulated as follows:
\begin{align}
    \min &\quad\displaystyle \sum_{a\in\calA^M}c^M_az_{ua} + \rho - \sum_{i\in\calI}\sum_{t\in\calT}\sum_{a\in\calA_{(i,t)}^{M,-}} \pi_iz_{ua} + \sum_{i\in\calI}\sum_{a=(i,t)\in\widetilde\calA^M_i}\xi_{it}z_{ua}
    \label{jo_jobssubproblem}\\
    \st&\quad\text{Constraints~\eqref{jobleaveonce}--\eqref{jobpartition}}\\
    &\quad z_{ua}=0,\ \forall i\in\calI,\ \forall t<T_i^S,\ \forall a\in\calA_{(i,t)}^{M,-}\\
    &\quad z_{ua}=0,\ \forall i\in\calI,\ \forall t>T_i^E,\ \forall a\in\calA_{(i,t)}^{M,+}\\
    &\quad \sum_{a\in\widetilde\calA^M_i}z_{ua}=\tau_i,\ \forall i\in\calI\\
    &\quad \text{$\bz$ binary}
\end{align}

This reformulation enables to cast the pricing problem as a shortest path problem with time windows in the physical network---as opposed to the time-expanded network. In this network, node 0 represents the starting depot, node $n+1$ represents the ending depot, and nodes $1,\cdots,n$ represent job locations. Note that the total cost of a path is decomposable into additive arc weights. This property is natural for the the travel costs, which can be added to each arc via the $d^M_{ij}$ parameters; for the fixed cost $\rho$, which can be added at the outset; and for the reward $\pi_i$, which can be added to any arc incoming into job $i\in\calI$. For the dual prices of the linking constraint $\xi_{it}$, however, there could be multiple periods when the job could be performed; for instance, a job with time window [4,8] and duration 3 could be performed between times 4 and 7, or between times 5 and 8. This situation can be captured in the physical network---without time expansion---by duplicating each arc into each node and assigning a corresponding cost of $\sum_{t\in\{t_0,\cdots,t_0+\tau_i-1\}}\xi_{it}$, where $t_0$ denotes the start time of the job for the arc under consideration. We shall see later on that we can restrict the search to those arcs that start the job as early as possible, thereby eliminating such arc duplicates.

\begin{figure}[h]
\centering
\includegraphics[scale=0.60]{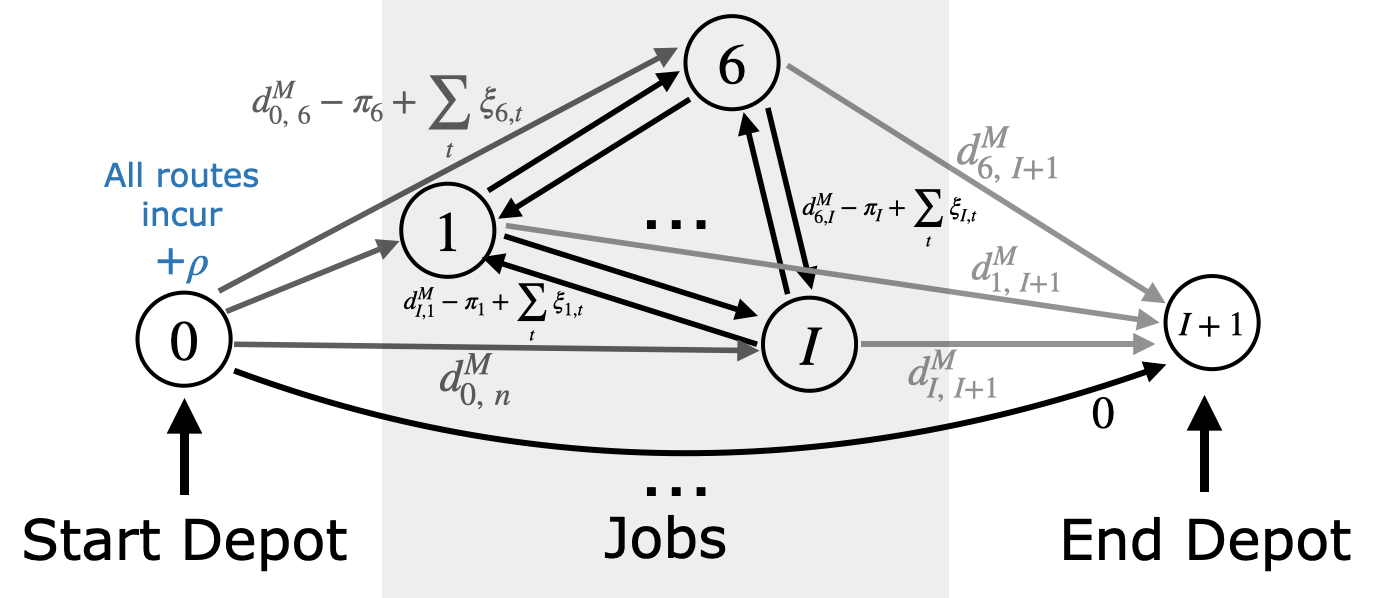}
\caption{Visualization of the graph structure to solve the pricing problem for mission vehicles. Each node represents a depot or a job location. The arc weights encode the travel time, the shadow cost of vehicle utilization, the shadow cost of coverage, and the dual reward of performing each job.}\label{fig:jobGraph}
\end{figure}

Figure~\ref{fig:jobGraph} shows the graph for the pricing problem in the physical network, along with the arc-based decomposition of the cost function.

We solve the pricing problem using a label-setting algorithm, described in Algorithm~\ref{alg:subproblem}. Throughout the algorithm, we maintain a list of paths with a bi-label structure to keep track of the time stamp (to enforce time windows) and the reduced cost (to minimize the objective function). The algorithm proceeds as follows:

\textit{State definition.} Let $(i^l,\mathbb{S}^l)$ denote a state, where $i^l$ tracks the ``current'' location, $T^l$ tracks the time stamp, and $\mathbb{S}^l\subseteq\{0,\cdots,n+1\}$ tracks the jobs that have been visited. For every state, we keep track of the time as $T(i^l)$ and of the reduced cost as $RC(i^l,\mathbb{S}^l)$.

\textit{Initial state:} $(i^0=0,\mathbb{S}^0= \emptyset)$, with $T(I^{0})={0}$, $RC(i^0,\mathbb{S}^0)=\rho$.

\textit{State transitions.} For any state $(i^l,\mathbb{S}^l)$, the state is potentially updated to $(i^{l'},\mathbb{S}^{l'})$, where $i^{l'}\in\{1,\cdots,n+1\}\setminus\mathbb{S}^l$ and $\mathbb{S}^{l'}=\mathbb{S}^l\cup\{i^{l'}\}$. We apply a feasibility check to apply time windows. Specifically, the transition is admissible if $\max(T(i^l)+TT^M(i^l,i^{l'}),T_{i^{l'}}^S)+\tau_{i^{l'}}-1\leq T_{i^{l'}}^E$, where $TT^M(i^l,i^{l'})$ denotes the minimum possible travel time between locations $i^l$ and $i^{l'}$.

\textit{Time label.} We update the time label by setting $T(i^{l'})=\max(T(i^l)+TT^M(i^l,i^{l'}),T_{i^{l'}}^S)$.

\textit{Reduced cost label.} We update the reduced cost label by setting:
$$\Pi(i^{l'}, \mathbb{S}^{l'})=\Pi(i^l,\mathbb{S}^l)+d^M_{i^l,i^{l'}}-\pi_{i^{l'}}+\sum_{t=\max(T(i^{l'})}^{T(i^{l'})+\tau_{i^{l'}}-1}\xi_{it}$$

\textit{Dominance rule.} We say that state $l^1$ dominates state $l^2$ if
\begin{enumerate*}[label=(\roman*)]
\item	$\mathbb{S}^{l^1}\subseteq \mathbb{S}^{l^2}$,
\item	$T(i^{l^1})\le T(i^{l^2})$, and
\item	$RC(i^{l^1}, \mathbb{S}^{l^1}))\le RCi(i^{l^2}, \mathbb{S}^{l^2})$.
\end{enumerate*}
Any dominated state is immediately pruned.

Upon termination, we return all non-dominated states. Pseudocode is provided in Appendix~\ref{appendix:spa}.

By design, the label-setting algorithm returns the shortest paths from the origin depot to all other nodes. If the smallest reduced cost is non-negative, the algorithm provides a certificate that the dual constraint of the set-partitioning formulation is satisfied for all routes---even omitted ones. Otherwise, the paths with negative reduced costs get added to $\calQ^0$ and the algorithm proceeds.

Note that the label-setting algorithm leverages the jobs' time windows to prune branches in the shortest path tree. However, some of the time windows may be loose enough so the mission vehicles still have flexibility to perform the same job at multiple times---even for a given sequence of jobs. For instance, consider a case where a route has three jobs to be done within time windows [1, 4], [5, 11], and [10, 13]; assume that all travel distances are equal to 1 and all work loads are equal to 3. Then, the second job can be done from time 5 to time 8, or from time 6 to time 9. To accelerate the algorithm, we restrict the search to those paths that perform the jobs as early as possible for a given sequence of jobs in the path. In the above example, the path performing the second job from [5, 8] would be considered, but the path performing the job at [6, 9] would not. Mathematically, we would start each job at $\max(T(i^l)+TT^M(i^l,i^{l'}),T_{i^{l'}}^S)$.

We propose an acceleration technique that processes all jobs as early as possible. Specifically, in the label-setting algorithm, we only consider those state transitions that  visit node $i^{l'}$ at time $\max(T(i^l)+TT^M(i^l,i^{l'}),T_{i^{l'}}^S)$. This technique, referred to as \textit{early job completion}, reduces the dimensionality of the problem and therefore accelerates computational times. In the vast majority of cases, the early job completion technique comes at no loss of optimality---although there may exist pathological cases in which it only approximates the original problem, in which case it becomes a heuristic.

\begin{proposition}\label{prop:earliestWorks}
    Early job completion is exact for the majority of cases, except special instances described below.
\end{proposition}

\begin{figure}[h]
\centering
\includegraphics[scale=0.50]{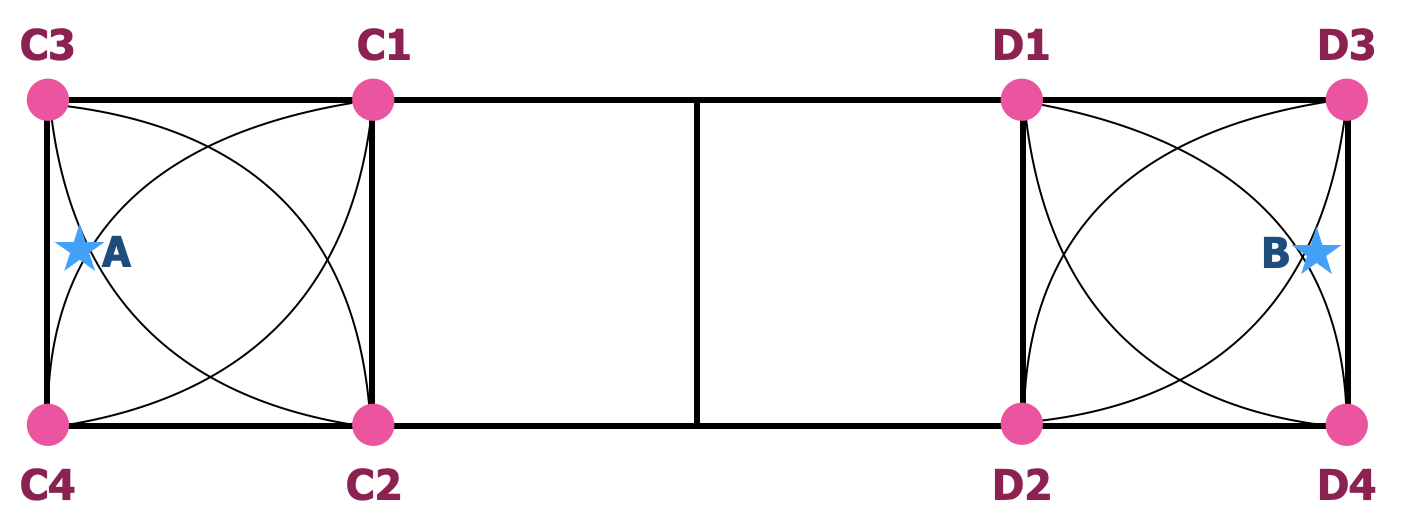}

\caption{A diagram of a pathological scenario, in which the job travel distance between $A$ and $B$ (blue stars) is one discrete time unit longer than the coverage travel distance between $C3$/$C4$ and $D3$/$D4$.}\label{fig:pathological}
\end{figure}

Specifically, an example of a case in which the early job completion technique fails to return the optimal solution is shown in Figure~\ref{fig:pathological}. In that case, the early job completion technique would invalidate a solution and force a new mission vehicle to be deployed to a certain job, whereas simply rescheduling that job to a later time would have solved the problem with the same objective as the would-be early job completion solution.

Suppose the fineness of the coverage mesh equalled the radius of coverage, and there are two jobs, A and B, just out of reach of coverage spots C1/C2 and D1/D2, respectively, forcing usage of coverage at C3/C4 and D3/D4. Also assume B has a large time window. (For ease, C refers to one of C3/C4 and D, one of D3/D4.) The distance between the viable $C$ and $D$ coverage spots is only barely greater than the distance between $A$ and $B$. Suppose the mesh has fineness 10, and the radius of coverage is 10. The distance between coverage is 40, while the distance between jobs is 37.33. Recall that travel time is given by a ceiling function, as it is discretized. In most cases, the speed will cause the time taken for jobs and coverage travels to be the same. But if the speed is a number such that the time difference between spots and jobs is nonzero, such as 19, then the coverage travel time is 3, while the jobs travel time is 2. Suppose further that job A finishes at time 4, and job B takes 5 time to finish, with time window [6, 19]. Then by early job completion, the mission vehicle rushes from A to B, reaching it at time 6, working from 6 to 11. But the emitting vehicle, which will be at C, needs 3 time to get to a valid D, so it can only start covering at time 7. Thus, this emitting vehicle will not be able to continue to cover the mission vehicle. A new emitting vehicle must be sent, increasing the distance cost beyond the optimum. Had the mission vehicle waited until 7 to reach B, it would have accompanied the mission vehicle and successfully completed the coverage. Even if the mesh were made finer, it is not guaranteed that these pathological combinations of vehicle speed and coverage spot location would disappear. One can imagine a scenario in which the optimal emitting and mission vehicle placement it the same as in Figure~\ref{fig:pathological}.

\subsection{Pricing problem: seeking paths for emitting vehicles}
\label{subsec:emitting}

The pricing problem for emitting vehicles seeks the path of minimal reduced cost, as follows:
\begin{align}
    \min &\quad\displaystyle \sum_{a\in\calA^E}c^E_ax_{va} + \beta - \sum_{j\in\calJ}\sum_{t\in\calT}\sum_{a\in\calA_{(j,t)}^{M,-}}\left(\sum_{i\in\calI:j\in\calC_i} \xi_{it}\right)x_{va}
    \label{jo_jobssubproblem}\\
    \st&\quad\text{Constraints~\eqref{coverageleaveonce}--\eqref{conserveflowcpv}; $\bx$ binary}
\end{align}

Again, this problem features a prize-collecting shortest path structure---without time windows. Each location $j\in\calJ$ is associated with a reward stemming from the shadow price of the linking constraint associated with \emph{all} job locations within the coverage radius. This expression enables the computation of the reward for each coverage location $j$, as opposed to each job location $i$. Each arc comes at a traveling cost $c_a^E$, which decomposes the total cost of the path $C^p$. Finally, the fixed cost $\beta$ captures the shadow cost of using one emitting vehicle.

This pricing problem can also be cast as a shortest path problem, and we solve it with a similar label-setting algorithm as shown in Section~\ref{subsec:mission}. One important difference, however, is that this shortest path formulation needs to be carried out in a time-expanded network. Indeed, the problem does not exhibit restrictions on the time during which an emitting vehicle stays at a coverage location. Therefore, the problem jointly optimizes the physical route of an emitting vehicle and the time at which the vehicle departs from each location and arrives at each subsequent location. Moreover, the number of candidate stopping locations for emitting vehicles is very large---as long as spatial discretization is sufficiently granular. Therefore, the time-space network can grow very large, which severely hinders the tractability of the pricing problem.

In response, we develop three acceleration strategies to speed up the search over the shortest path tree, by leveraging the linking constraint between mission vehicles and emitting vehicles. First, we restrict the search to the pairs of coverage locations and times that can cover at least one job---we refer to this technique as \textit{input-based pruning}. Specifically, an emitting vehicle is only allowed to enter a coverage spot $j\in\calJ$ at time $t\in\calT$ if at least one job within the coverage radius starts at some time $t$ (i.e. the first of the work is being done at time $t$), and no other jobs are in progress at time $t$ (i.e. its time being worked, as determined by the job routing, cannot have  unless it is the start). Similarly, a emitting vehicle can only leave a coverage spot at time $t$ if at least one job within the coverage radius ends at time $t$ (i.e. the last of the work is being done at time $t$), and no other jobs are in progress at time $t$ (i.e. its time being worked, as determined by the job routing, cannot include $t$ unless it is the end). In the case that a mission plan is not available, the job time windows will be used instead of mission vehicle work windows.

\begin{proposition}\label{prop:inputBased}
Input-based pruning does not imply any loss of optimality.
\end{proposition}

This input-based pruning technique reduces the number of possible times for entry and exit of a spot to a few distinct possibilities governed by when jobs are being performed. Further window pruning can be performed by further restricting the times of entry and exit in a coverage location to those times when a job is being performed in the coverage area, based on the incumbent RMP solution---we refer to this technique as \textit{primal-based window pruning}. Mathematically, consider a coverage spot $j$ which potentially covers jobs $\calI_j = \{ I_{j1}, I_{j2}, \dots, I_{jm_j} \}$ with associated time windows $[T_{j1}^S, T_{j1}^E], \dots, [T_{jm_j}^S, T_{jm_j}^E]$. Then, the set of potential entry times is comprised of all values $T_{ji}^S$ for which $\forall (i \neq ji) \in \calI_j, \ T_{i} \notin (T_{i}^S, T_i^E]$. In other words, for any entry time, no job within the range of this coverage spot is concurrently performed. Similarly, the set of potential exit times equals the job end times for which no concurrent jobs are being performed. For example, if a coverage spot has jobs with work windows $[4, 7], \ [6, 9], \ [10, 12]$, then the valid entry times are $\{ 4, 10\}$, while the valid exit times are $ \{ 9, 12 \}$.

Unlike input-based pruning, the window pruning technique is a heuristic. Indeed, it restricts the search in the pricing problem for emitting vehicles to a neighborhood defined by the incumbent primal RMP solution for mission vehicles. However, the incumbent RMP solution may not be optimal, and there may beneficial time-location pairs that are not visited under this heuristic. For this reason, we apply window pruning at the beginning of the algorithm to accelerate it, and then resort to the exact pricing problem algorithm to guarantee optimality upon termination.

Finally, we propose another acceleration strategy, referred to as \textit{dual-based pruning}. This technique is motivated by the observation that few of the linking constraints are binding at the optimum, and, accordingly, few of the $\xi_{it}$ variables are greater than zero. We leverage this observation to prove that the solution of the pricing problem for emitting vehicles will never visit a location $j\in\calJ$ at time $t\in\calT$ if all jobs $i\in\calI$ in its coverage area (i.e., $j\in\calC_i$) are such that $\xi_{it}=0$. Indeed, all other terms in the reduced cost expression are positive, which creates no incentive to visit any location if the ``reward'' is zero. This enables us to restrict the search to intervals with at least one positive $\xi_{it}$ in the vicinity of location $j$, thereby considerably enhancing the scalability of the algorithm.

Specifically, the dual-based pruning technique proceeds by identifying all beneficial job-time pairs, stored in a set $\calB_j=\{(i,t)\in\calI\times\calT:\xi_{it} > 0\, \ i \ \text{covered by} \ j\}$. In other words, $\calB_j$ is associated with a coverage spot, which covers many jobs $i\in\calI$. For any coverage location $j \in \calJ$, we collect, for each job $i \in \calI$, all time intervals $\calT_i$ for which there exists at least one time $t \in T_i \in \calT_i$ with $\xi_{it} > 0$ ($T_i$ is one time interval within $\calT_i$). Let $\mathcal{X}$ denote the possible sets of $\xi_{it} > 0$ values included in any $T_i$. We trim down $\calT_i$ so that for each set $X \in \mathcal{X}$, only the unique timewise-shortest interval $T_i$ containing exactly the $\xi_{it}$ in $X$ is included. Let $\calT^*_j$ be the union of all $\calT_i$ for all $i$ within coverage range of $j$. The label-setting algorithm merely iterates over spots $j$ and then $\calT^*_j$, rather than spots $j$ and then possible arrival and departure times.

\begin{proposition}\label{prop:dualBased}
Dual-based pruning does not imply any loss of optimality.
\end{proposition}

\subsection{Restoring integrality: a stem-and-blender approach}

At this point, we have presented our algorithm to solve the linear relaxation of the set-partitioning formulation. However, this does not preclude the possibility of fractional solutions at the optimum. One possibility to eliminate them would be to proceed via branch-and-price, albeit at significant costs in terms of computational times. Instead, we propose a perturbation-scheme to increase the pool of path-based solutions in the restricted master problem, and solve a final binary optimization problem over this augmented set of candidates.

The motivation for this approach lies in the regular patterns observed in the fractional solutions. Let us denote by $\calQ^*$ and $\calP^*$ the sets of paths for mission vehicles and emitting vehicles at the optimum of the set-partitioning relaxations, respectively. For the sake of the example, assume that a path-based solution involves two routes $q_1$ and $q_2$ with weights $\overline{z}_{q_1}=\overline{z}_{q_2}=0.5$. Then, the two paths need to share a significant fraction of the route, with some deviations at some point in time. Accordingly, we divide each path into a shared ``stem'' and a non-shared ``tail''. We then partition the set $\calQ$ (resp. $\calP$) into ``stem sets" $S_1, S_2, \dots, S_{|S|}$ where all routes in $S_s$ share the same first job location (resp. coverage location). the stem of each group is then defined as the longest set of job locations (resp. coverage locations) that are shared across all paths in $S_s$; the blender lists all subsequent locations visited by at least one path.

\textit{Example.} Suppose $\calQ^*$ is given by the following paths, with weights $\overline{z}_q$ on the left:

0.25: 01 63 64 65 04 05 06 15

0.25: 01 63 64 65 10 11 12 13 14 15

0.50: 01 63 03 66 67 68 69 70 71

0.25: 16 73 02 03 66 67 68

0.25: 16 73 02 03 66 67 68 69 70 71 

0.50: 16 73 02 04 05 06

0.25: 21 07 08 09 65 04 05 06 69

0.25: 21 07 08 09 65 10 11 12 13 14 70 71 

0.50: 21 07 08 64 09 10 11 12 13 14 15

In this example, the stems are: [01 63], [16 73 02], and [21 07 08]. Note that [01 63 64 65] is not a stem given how stems are defined above: fractional routes are partitioned by first job, and then any stem must encapsulate all routes of a particular partition. The blender associated with the stem [21 07 08] consists of elements: 04, 05, 06, 09, 10, 11, 12, 13, 14, 15, 64, 65, 69, 70, 71.

The stem-and-blender approach then considers all the nodes in each stem (corresponding to each stem set $S_s$), and runs the label-setting algorithm (Sections~\ref{subsec:mission}--\ref{subsec:emitting}) to generate all feasible paths starting with the stem and visiting elements of the blender thereafter. In other words, this approach restricts the search to the patterns encountered across all paths from the linear relaxation solution, but introduces perturbations from the label-setting algorithm to increase diversity in the path sets. The resulting paths are then added to $\calQ^*$ and $\calP^*$, and we run the restricted master problem one last time while enforcing integrality.

\subsection{Full DCG algorithm}\label{sec:fullAlgorithm}

The full DCG algorithm is detailed in Algorithm~\ref{alg:dcg}. We initialize feasible path sets $\calQ^0$ and $\calP^0$ for mission vehicles and emitting vehicles, respectively. The algorithm iterates between the RMP and both pricing problems. At each iteration, the algorithm adds up to one path for mission vehicles to $\calQ^0$ and up to one path for emitting vehicles to $\calP^0$. The algorithm stops when there is no path with strictly negative reduced cost from both pricing problems. Upon termination, we call the stem-and-blender procedure to augment the sets $\calQ^0$ and $\calP^0$, and solve the RMP one last time with integrality constrains. The algorithm yields a feasible solution to the set-partitioning formulation.

\begin{algorithm}[h!]
\floatname{algorithm}{Algorithm}
\caption{Double Column Generation (DCG).}
\renewcommand{\arraystretch}{1.0}\small
\label{alg:dcg}
\begin{algorithmic}[1]
\Require{Initial mission routes $\calQ^0$, initial emitting routes $\calP^0$, possibly through warm starts.}
\State Perform initial calculation of parameters: $C^q, \ C^p, \ \calQ_i, \ \delta_{it}^q, \ g_{jt}^p$.
\While{\text{true}} \Comment{Stop if neither subproblem returns new route}
\State Solve Restricted Master Problem (RMP); extract dual variables $\pi_i, \ \rho, \ \beta, \ \xi_{it}$
\State Solve shortest paths pricing problem for mission vehicles
\State Solve shortest paths pricing problem for emitting vehicles
\If{\text{at least one new route generated}}
\State Recalculate $C^q, \ C^p, \ \calQ_i, \ \delta_{it}^q, \ g_{jt}^p$
\Else
\State Leave loop
\EndIf		
\EndWhile
\State Perform stem-and-blender on $\calQ^0$ and $\calP^0$
\State Solve RMP one last time with integrality constraints
\end{algorithmic}
\end{algorithm}

We apply a warm start technique to generate initial plan sets $\calQ^0$ and $\calP^0$. Specifically, we perform a mission-only column generation algorithm, and set all considered routes in that algorithm to $\calQ^0$; then, based on the mission-only solution, we perform an emitting-only column generation algorithm, and set all considered routes in that algorithm to $\calP^0$.

Finally, we apply the four acceleration strategies described in this section: \textit{early job completion}, \textit{input-based pruning}, \textit{primal-based pruning}, and \textit{dual-based pruning}. Note that all but primal-based pruning are exact acceleration methods, which do not imply any loss of optimality; we therefore apply them throughout the algorithm. We apply \textit{primal-based pruning} in early iterations, and then de-activate it to generate a certificate of optimality upon termination, after which we call the stem-and-blender procedure to augment the sets $\calQ^0$ and $\calP^0$, and solve the RMP one last time with integrality constrains. 

\section{Computational results}
\label{sec:results}

We implement the proposed model and the DCG algorithm to realistic problem inputs, developed in collaboration with Singapore's Defense Science and Technology Agency (DSTA). The instances are designed to replicate mission-critical environments, with clustered jobs performed in remote locations where broadband internet is otherwise unavailable. Specifically, we define an experimental setup based on four parameters:
\begin{itemize}
    \item[--] number of jobs to be performed within the time horizon
    \item[--] number of job clusters: the more clusters, the fewer jobs each cluster contains.
    \item[--] cluster radius: the larger the radius, the smaller the incidence of clusters in the problem. 
    \item[--] radius of coverage: range of WiFi coverage from each emitting vehicle.
\end{itemize}

We consider a rectangular area in which vehicles can evolve freely. For each cluster, we uniformly sample the centroid in the overall area, and uniformly sample job locations within the Euclidean ball defined by the centroid and the radius. In particular, the cluster radius and the coverage radius determine the interdependencies between the operations of emitting vehicles and mission vehicles.

\begin{figure}[h]

\centering
\includegraphics[scale=0.4]{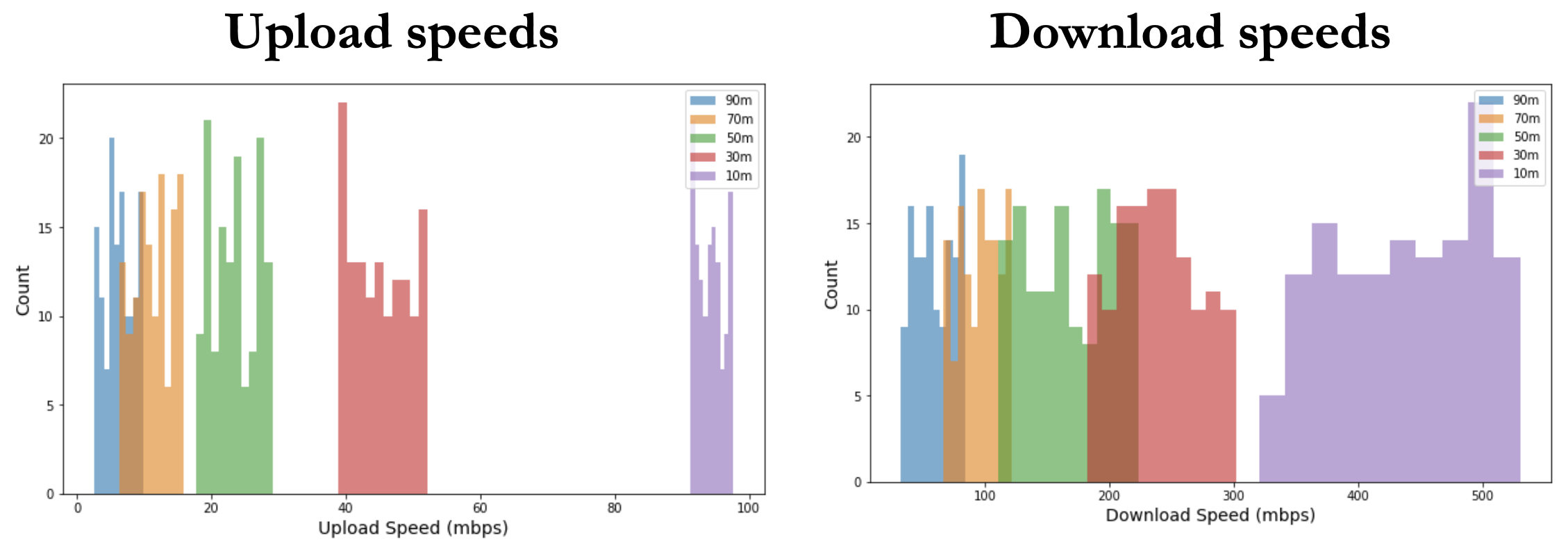}
\caption{Upload (left) and download (speeds) for 5G with the specified radius from the mast. Over 100 samples were collected for each of the ten configurations.}
\label{fig:speed}
\end{figure}

We gathered real-world data on WiFi coverage from emitting vehicles, in collaboration with the DSTA. Specifically, a 5G mast was deployed in several locations, and a car equipped with a mobile device was traveling around the emitter. We collected data that included the distance from the 5G mast and recorded upload and download speeds. Figure~\ref{fig:speed} reports the histogram of the upload speeds as a function of distance---note that, with 5G, download speeds are over five times faster than upload speeds. Based on these data, we use a coverage radius of 50 in our experiments, with coverage mesh size 50 as well.

All models are solved with Gurobi v8.1, using the optimization package JuMP in Julia \citep{dunning2017jump}. All runs were performed on a laptop equipped with a 4-core Intel i5 CPU (2 GHz) and 16 GB RAM, with a two-day limit. All instances and code are available online.\footnote{\hyperlink{https://github.com/DA72p/ict-opt/blob/main/full-solution-notebook.ipynb}{https://github.com/DA72p/ict-opt/blob/main/full-solution-notebook.ipynb}}

\subsection{Scalability of algorithm}

We first assess the effectiveness of the column generation algorithm for solving the problem of mission vehicles. In this experiment, we merely optimize the operations of mission vehicles following a vehicle routing structure with time windows without consideration for network coverage. Table~\ref{tab:benefitsOfCG} reports the average computational times over 50 random trials, as a function of the the number of jobs. In these experiments, we consider five dense clusters---corresponding to a small cluster radius.

\begin{table}[h!]
\centering
\caption{Computational times (in seconds) of each method for optimizing operations of mission vehicles.}\label{tab:benefitsOfCG}
\begin{tabular}{ccccccc}
\toprule[1pt] \vspace{1pt}
&&&\multicolumn{3}{c}{Column generation}\\\cmidrule{4-6}
Num. Jobs & Explicit     & Set partitioning     & Naïve     & with domination & with Proposition~\ref{prop:workwhenstay} \\ \hline
15        & 13.59  & 3.56   & 1.55   & 1.57   & 0.15       \\ %\hline
20        & 26.02  & 6.07   & 2.42   & 2.53   & 0.48       \\ %\hline
25        & 128.93 & 14.56  & 7.81   & 4.05   & 1.05       \\ %\hline
30        & 540.39 & 30.05  & 22.70  & 12.16  & 2.13       \\ %\hline
35        & ---       & 224.39 & 145.49 & 23.45  & 5.76       \\ %\hline
40        & ---       &  ---      & ---       & 35.56  & 8.90       \\ %\hline
45        & ---       &  ---      & ---       & 77.84  & 14.56      \\ %\hline
50        & ---       &  ---      & ---       & 256.49 & 44.69      \\ %\hline
55        & ---       &  ---      & ---       &  ---      & 62.79      \\ %\hline
60        & ---       &  ---      & ---       &  ---      & 151.46     \\ %\hline
65        & ---       &  ---      & ---       &  ---      & 365.89     \\ %\hline
\bottomrule[1pt]
\end{tabular}
\end{table}

Note, first, the benefits of column generation over the explicit formulation and the full set-partitioning formulation. Whereas off-the-shelf algorithms can handle small-scale instances with up to 30--35 jobs, they time out  as the problem becomes larger. In comparison, column generation can handle medium-scale instances. Yet, column generation by itself is insufficient to scale to large-sized instances, thus motivating our acceleration strategies. First, using domination criteria for tree pruning in the label-setting algorithm cuts down computational times by nearly one order of magnitude with 35 jobs and scales to instances with 50 jobs. Then, further gains are achieved by leveraging the fact that mission vehicles do not have to stay at a job location when not performing work (Proposition~\ref{prop:workwhenstay}). Again, this result enables speedups by nearly one order of magnitude and further increases the scalability of the solution method to large instances with 65 jobs. These results are robust across parameter settings, captured via different numbers of clusters and cluster radii.

Next, we implement the DCG algorithm to solve the integrated problem---jointly optimizing the operations of mission vehicles and emitting vehicles subject to linking constraints ensuring network coverage when mission vehicles are at the job locations. Table~\ref{tab:benefitsOfWS} reports average computational times over 50 random seeds for each number of jobs, in the same setting as in Table~\ref{tab:benefitsOfCG}. We first consider a baseline implementation of the algorithm, which includes the label-setting algorithm with the domination criteria along with a restriction to visit jobs only when performing work, following Proposition~\ref{prop:workwhenstay}. As such, the operations of the mission vehicles under the baseline implementation in Table~\ref{tab:benefitsOfWS} corresponds to the rightmost column in Table~\ref{tab:benefitsOfCG}. We then add the acceleration strategies presented in Section~\ref{sec:algorithm} and the warm start (we implement partial warm starts for mission vehicles only, and full warm starts for both types of vehicles).

\begin{table}[h!]
\centering
\caption{Computational times (in seconds) of the DCG algorithm with acceleration strategies.}\label{tab:benefitsOfWS}
\begin{tabular}{ccccc}
\toprule[1pt] \vspace{1pt}
Num. Jobs & Baseline   & with accelerations & with partial warm start & with warm start \\ \hline
15        & 1.88  & 1.19   & 0.78   & 2.12    \\ %\hline
20        & 16.83  & 3.15   & 1.53   & 2.62    \\ %\hline
25        & 95.68 & 18.42  & 2.36   & 3.23    \\ %\hline
30        & --- & 23.40  & 6.05   & 13.34   \\ %\hline
35        & --- & 56.71  & 12.57  & 19.34   \\ %\hline
40        & ---      & 280.69 & 32.09  & 48.07   \\ %\hline
45        & ---      &  ---      & 74.49  & 78.46   \\ %\hline
50        & ---      &  ---      & 150.67 & 133.33  \\ %\hline
55        & ---      &  ---      & 295.09 & 231.56  \\ %\hline
\bottomrule[1pt]
\end{tabular}
\end{table}

The DCG algorithm scales to medium-scale instances with dozens of jobs, but nonetheless fails to solve larger instances within the time limit. This limitation stems from the complexity of linking the operations of both types of vehicles in the master problem---in practice, it may take time for the dual values $\xi_{it}$ to converge to their optimal values, which can lead to a large number of iterations with weak coupling between mission and emitting vehicles. The acceleration strategies strongly speed up the algorithm and enable the solving of larger cases. Decomposing the acceleration techniques one by one yields input-based pruning as by far the most vital acceleration technique, as the algorithm could not even run for more than 20 jobs if this technique was not included. Dual-based pruning and the heuristic of primal-based pruning were helpful, but not as much as input-based pruning. Then, the warm starts are also instrumental to enhance the scalability of the algorithm. In small-scale instances, the full warm start can actually slow down the algorithm by increasing the computational requirements of the master problem with limited benefits in terms of a reduced number of iterations. As the problem size increases, however, warm starting the algorithm with close-to-optimal routes ends up providing computational benefits. With the warm start, the algorithm generates stronger primal solutions in early iterations and, most importantly, stronger approximations of the optimal dual variables $\xi_{it}$ based on the initial set of paths for mission and emitting vehicles. Ultimately, this stronger initialization enables faster convergence in the DCG algorithm.

\begin{figure}[h]
\centering
\includegraphics[scale=0.35]{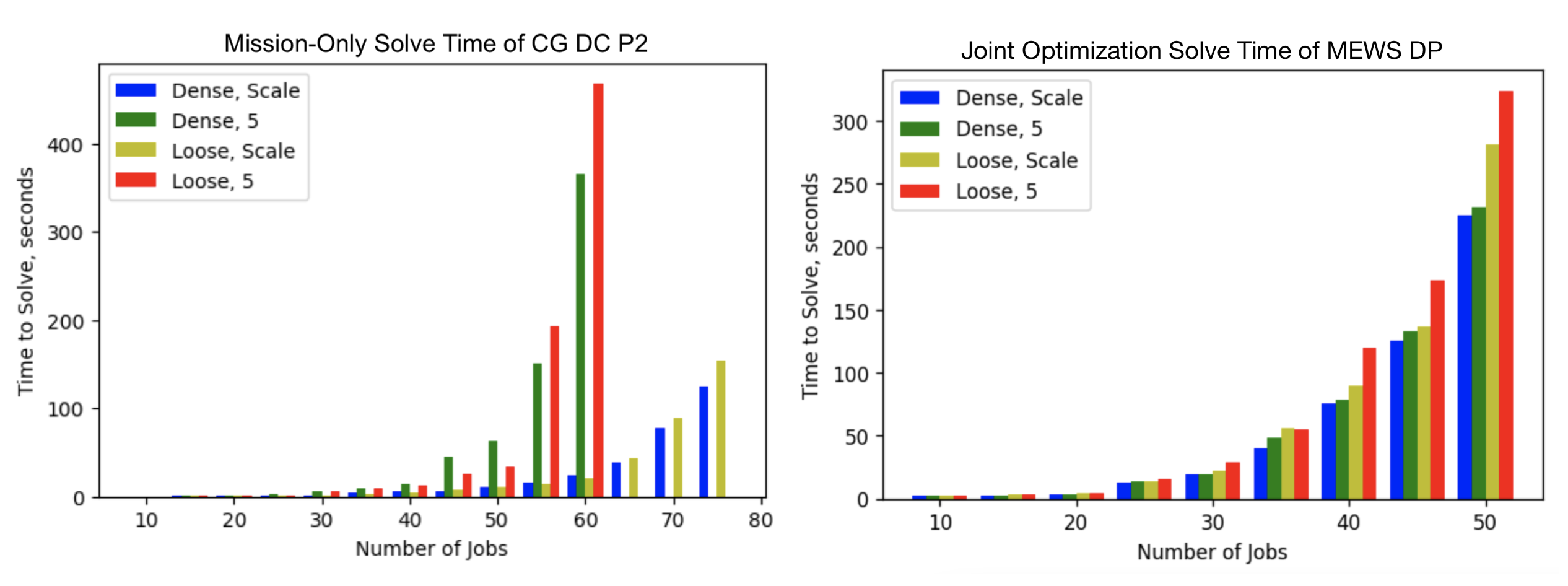}

\caption{(left) Graph of mission-only time. After 60 jobs, the 5-cluster cases take too long to compute and are not included in the graph. (right) Graph of joint time.}\label{fig:coverageScenarios}
\end{figure}

Figure~\ref{fig:coverageScenarios} reports the computational times of the best algorithms for the mission-only problem (left) and the full integrated problem (right), as a function of the size of each cluster (parameterized by the number of clusters) and the density of each cluster (parameterized by the radius of each cluster). Cluster radius is varied 20 (dense cluster) or 60 (loose cluster). The number of clusters either scales with the number of jobs (five jobs per cluster), or stays constant at 5.

As expected, the problem becomes increasingly complex as the number of jobs increases. In addition, the figure underscores the importance of the geometric setting on problem scalability. First, dense clusters render the problems more computationally tractable than loose clusters. This is intuitive, in that larger clusters lead to more complex routing problems, both for mission vehicles and for emitting vehicles. This effect is particularly important in the integrated problem due to the complex coordination requirements between both types of vehicles across wider geographic locales. In the mission-only problem, another critical driver of computational complexity is the number of clusters: the problem becomes much more complex as the number of jobs per cluster increases---that is, when the number of clusters remains constant at 5---and remains comparatively easier with more clusters and therefore fewer jobs per cluster. Indeed, larger clusters involve more complex routing problems in a vehicle routing structure with time windows. In the integrated problem, however, the number of clusters has a smaller effect. This can be explained by the fact that the operations of emitting vehicles are easier with fewer clusters---emitting vehicles can seek central locations in each cluster to support mission vehicles. With more clusters, however, the joint problem involves more complex coordination challenges between mission vehicles and emitting vehicles across clusters.

Finally, we evaluate the impact of the stem-and-blender approach to eliminate fractional routes in Table~\ref{tab:sab}, using the same setting as in Table~\ref{tab:benefitsOfWS}. As expected, the DCG algorithm can leave an optimality gap if we seek a solution from the variables generated during the column generation process---in particular, a relatively high fraction of routes take fractional values for emitting vehicles. However, the stem-and-blender approach eliminates \textit{all} fractional routes, therefore guaranteeing the optimality of the solution in the discrete optimization problem under consideration. These results underscore the impact of the stem-and-blender heuristic to solve discrete optimization problems with an exponential number of variables without resorting to more complex and potentially more time-consuming procedures such as branch-and-price.

\begin{table}[h!]
\centering
\caption{Average performance of the algorithm before and after applying stem-and-blender approach.}\label{tab:sab}
\begin{tabular}{cccccc}
\toprule[1pt]
Num. Jobs & S\&B Time & \# (\%) FMR & \# (\%) FER & \# FR post-S\&B & Opt. Decr.\\ \hline
45        & 10.58  & 1.57 (17.2)   & 3.33 (41.2)    & 0 & -0.16\% \\ %\hline
50        & 29.33  & 9.88 (65.5)   & 5.14 (55.6)    & 0 & -0.63\% \\ %\hline
\bottomrule[1pt]
\end{tabular}

\begin{tablenotes}\linespread{1}
\vspace{-6pt}
\item ``S\&B Time'': How long it takes to perform stem and blender, seconds\vspace{-3pt}
\item ``\# (\%) FMR": Number and percentage of fractional mission routes pre-stem \& blender\vspace{-3pt}
\item ``\# (\%) FER": Number and percentage of fractional emitting routes pre-stem \& blender\vspace{-3pt}
\item ``\# FR post-S\&B": Number of fractional mission and emitting routes after stem and blender\vspace{-3pt}
\item ``Opt. Decr.": Decrease in optimal objective after applying stem and blender\vspace{-3pt}
\end{tablenotes}
\end{table}

\subsection{Benefits of Optimization}

We demonstrate the benefits of optimization over baseline solutions that could be obtained in the absence of the methodology developed in this paper. Specifically, we consider three benchmarks:
\begin{itemize}
    \item[--] \textit{Greedy algorithm.} Each mission vehicle is sent out to visit the closest available job at any time, until time or jobs run out. We repeat this procedure until all jobs are covered by one vehicle. We then assign one emitting vehicle to each job route, and each emitting vehicle will travel to the closest location that covers the next job. This benchmark involves no optimization and provides a practical benchmark that could be easily applied in practice.
    \item[--] \textit{Mission-first optimization with one following emitting vehicle (``Emitting-follow").} This benchmark first optimizes the mission-only problem, by solving a vehicle routing problem with time windows without consideration for network coverage. In a second step, the benchmark assigns one emitting vehicle per mission vehicle under the same greedy strategy as above---namely, each emitting vehicle follows the mission vehicle perfectly while minimizing distance traveled. This benchmark corresponds to a heuristics based on a standard vehicle routing problem.
    \item[--] \textit{Mission-first, emitting-second optimization (``Mission-Emitting").} This benchmark applies the same first step as the previous one, by optimizing the operations of the mission vehicles without network coverage considerations. In a second step, the benchmark optimizes the operations of emitting vehicles, under the constraint that each mission vehicle must be covered when at a job location. This benchmark corresponds to an decomposition heuristics that could be implemented by solving one vehicle routing problem at a time.
\end{itemize}

Table~\ref{tab:performanceOfAlgos} shows average results over 20 random seeds. These experiments consider a setting with 50 jobs, three values of the cluster radius, and a constant coverage radius. Specifically, the table reports the solution obtained with each method (``Solution''), the distance ratio between mission vehicles and emitting vehicles (``DR''), the number of mission vehicles and emitting vehicles (``MV'' and ``EV''), and the number of coverage locations visited by emitting vehicles (``CL'').

\begin{table}[h!]
\centering
\caption{Average performance metrics for the joint optimization solution and the benchmarks.}\label{tab:performanceOfAlgos}
\begin{tabular}{c ccccc ccccc}
\toprule[1pt] \vspace{1pt}
       & \multicolumn{5}{c}{Cluster radius 25} & \multicolumn{5}{c}{Cluster radius 50}  \\ \cmidrule(lr){2-6}\cmidrule(lr){7-11}
Algorithm         & \multicolumn{1}{c}{Solution}   & \multicolumn{1}{c}{DR}   & \multicolumn{1}{c}{\#MV}   & \multicolumn{1}{c}{\#EV}   & \#CL   & \multicolumn{1}{c}{Solution}   & \multicolumn{1}{c}{DR}   & \multicolumn{1}{c}{\#MV}   & \multicolumn{1}{c}{\#EV}   & \#CL   \\ \hline
Greedy     & \multicolumn{1}{c}{167} & \multicolumn{1}{c}{1.15} & \multicolumn{1}{c}{10.2} & \multicolumn{1}{c}{10.4} & 3.4 & \multicolumn{1}{c}{186} & \multicolumn{1}{c}{1.11} & \multicolumn{1}{c}{11.4} & \multicolumn{1}{c}{11.3} & 3.6 \\
Emitting-follow & \multicolumn{1}{c}{113} & \multicolumn{1}{c}{1.26} & \multicolumn{1}{c}{9.8}  & \multicolumn{1}{c}{9.8}  & 1.5 & \multicolumn{1}{c}{123} & \multicolumn{1}{c}{1.24} & \multicolumn{1}{c}{10.1} & \multicolumn{1}{c}{10.1} & 3.0 \\
Mission-Emitting   & \multicolumn{1}{c}{102} & \multicolumn{1}{c}{1.49} & \multicolumn{1}{c}{12.2} & \multicolumn{1}{c}{8.2}  & 1.2 & \multicolumn{1}{c}{109} & \multicolumn{1}{c}{1.59} & \multicolumn{1}{c}{11.8} & \multicolumn{1}{c}{8.2}  & 2.5 \\ 
Joint optimization      & \multicolumn{1}{c}{100}    & \multicolumn{1}{c}{1.64} & \multicolumn{1}{c}{9.9}  & \multicolumn{1}{c}{8.1}  & 1.2 & \multicolumn{1}{c}{100}    & \multicolumn{1}{c}{1.67} & \multicolumn{1}{c}{10.1} & \multicolumn{1}{c}{8.1}  & 2.4 \\
\bottomrule[1pt]
\end{tabular}
\begin{tablenotes}\linespread{1}
\vspace{-6pt}
\item ``Solution'': objective value of each solution, normalized with respect to the optimal solution\vspace{-3pt}
\item ``DR'': distance ratio between mission vehicles and emitting vehicles\vspace{-3pt}
\item ``MV'', ``EV'': number of mission vehicles and emitting vehicles, respectively\vspace{-3pt}
\item ``CL'': number of coverage locations visited by the emitting vehicles
\end{tablenotes}
\end{table}

Note, first, that the greedy algorithm induces significant optimality losses as compared to joint optimization, with 60--100\% cost increases. The optimization-based benchmarks perform comparatively much better, but nonetheless result in inferior solutions. In particular, the two-stage decomposition method results in a 2--9\% loss. In fact, the decomposition heuristic tends to perform worse as the number of jobs increases, indicative of harder problems with stronger interdependencies between mission vehicles and emitting vehicles. In comparison, the joint optimization solution consistently returns the optimal solution in instances with up to 50--60 jobs, in manageable computational times (as discussed in the previous section). As a result, the joint optimization methodology developed in this paper can provide significant practical benefits toward coordinating ICT operations and physical operations in emerging connected mission-critical environments.

These results further shed light into the drivers of the benefits of optimization. Specifically, the joint optimization solution achieves a higher distance ratio between mission vehicles and emitting vehicles---by up to 10\% as compared to the second-best solution. In other words, the joint optimization solution achieves coverage for all mission vehicles with a smaller distance traveled by emitting coverage, therefore achieving higher synergies between the two types of vehicles. Along the same lines, the joint optimization solution also leverages fewer vehicles and visits fewer coverage locations, again suggesting more efficient coordination between mission and emitting vehicles. Ultimately, the benefits of optimization mainly come from aggregating the operations of mission vehicles in time and space to cover them efficiently with a small number of vehicles at limited costs in terms of distance traveled---especially for emitting vehicles.

\begin{figure}[h]
\centering
\includegraphics[scale=0.55]{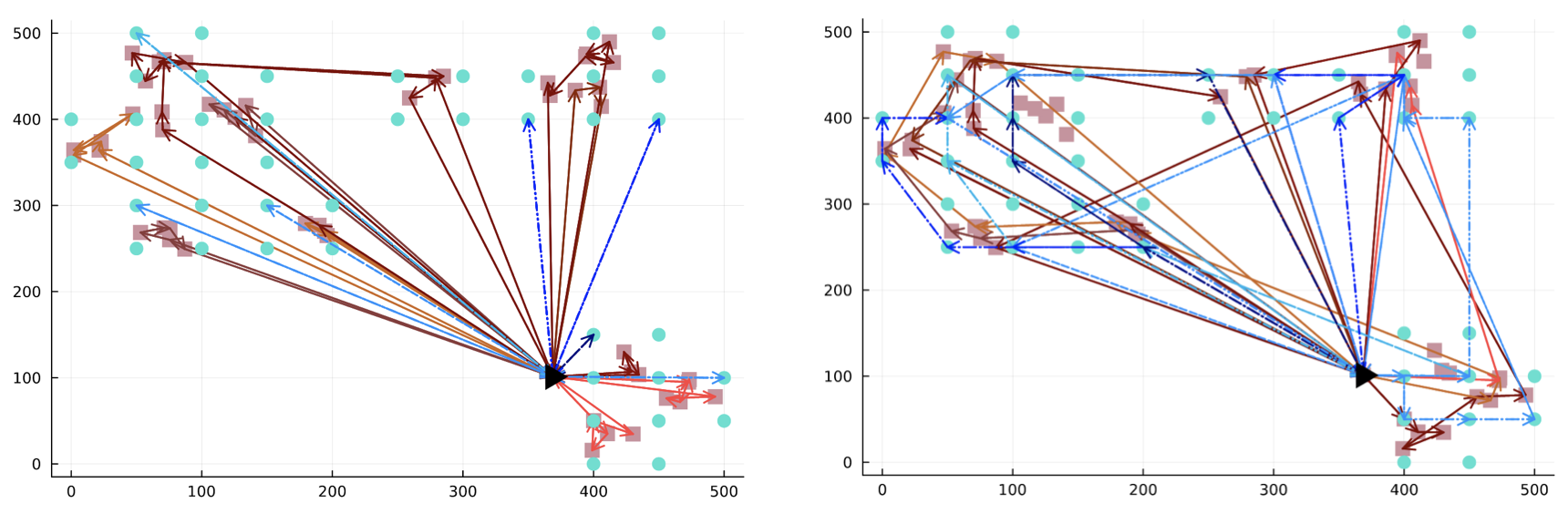}
\caption{Joint optimization solution (left) and greedy solution (right), for 50 jobs and cluster radius 25. Each plot is a physical 2d location map, with the y-axis representing latitude and the x-axis representing longitude. Aqua dots represent candidate coverage spots, while mauve squares represent job locations. The black triangle is the depot location. Red or brown routes are mission vehicle trajectories, while blue routes are emitting vehicle trajectories.}
\label{fig:map}
\end{figure}

Figure~\ref{fig:map} illustrates the joint optimization solution and the greedy benchmark. As expected, emitting vehicles travel much further distances in the greedy solution, leading to a solution deterioration of 65\% in that case. For these reasons, we also find that the benefits of optimization increase as the cluster radius decreases and as the coverage radius increases (i.e., jobs get more concentrated, thereby facilitating the operations of emitting vehicles).

\subsection{Benefits of Mobile Emitters}

Recall that mobile emitting vehicles (Figure~\ref{fig:5gtruck}) are highly expensive. This motivates the exploration of solutions with fewer emitting vehicles, potentially at a cost in terms of distance traveled. Through this analysis, we compare two optimization solutions, corresponding to two practical instances: (i) mobile emitters, using the joint routing optimization method developed in this paper; and (ii) fixed emitters, corresponding to an instance where emitting vehicles need to be located in the same place throughout the planning horizon. Specifically, the fixed emitter scenario is solved via a two-stage approach, which first minimizes the number of emitting vehicles to cover all the jobs, and then optimizes the routes of mission vehicles subject to the resulting coverage.

\begin{figure}[h]
\centering
\includegraphics[scale=0.5]{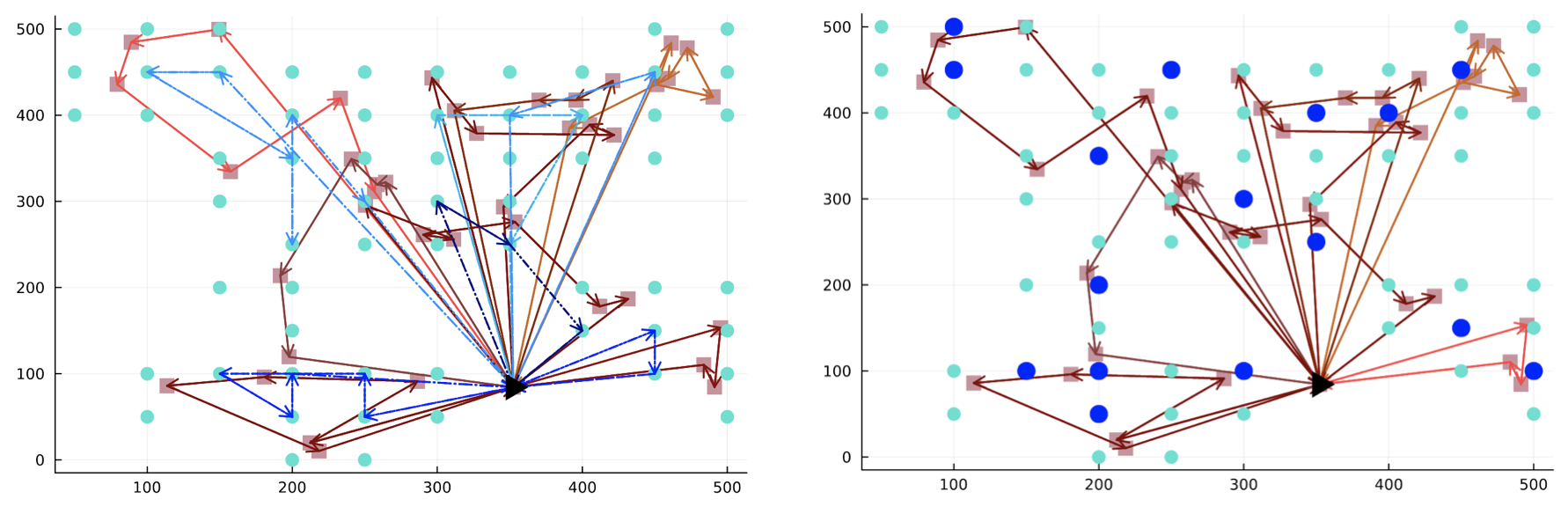}

\caption{Solution with mobile emitting vehicles (left) and stationary emitting vehicles (right). The black triangle is the depot, mauve squares are job locations, aqua circles are coverage spots, and shades of brown are job routes. There are 40 jobs, coverage radius is 50, and cluster radius is 100. Left: There are 6 coverage routes (various shades of blue). Right: The 16 blue dots denote stationary emitting vehicles.}
\label{fig:map_fixed}
\end{figure}

Table~\ref{tab:necessaryEmitting} reports the minimum number of emitting vehicles required to cover all missions, averaged over 50 seeds for each number of jobs (40 and 50) and each cluster radius (25, 50, and 100). Note, first and foremost, that mobile emitters result in a significantly smaller fleet than stationary emitters. The difference becomes larger as the cluster radius increases. Indeed, when jobs are physically disperse, synergies across job locations become weaker and it is critical for an emitting vehicle to move from one location to another to cover multiple jobs. In the most disperse cases, mobile emitting vehicles can reduce the fleet size by a factor of 2--3. This is illustrated in Figure~\ref{fig:map_fixed} in an instance with a large cluster radius of 100. note that all six mobile emitting vehicles move to at least two spots. Ultimately, the methodology developed in this paper highlights benefits of mobile emitting technologies by providing added flexibility to the routing problem and enabling operators to perform more mission-critical jobs with smaller fleets of emitting vehicles.

\begin{table}[h!]
\centering
\caption{Average number of emitting vehicles required, as a function of number of jobs and cluster radius.}\label{tab:necessaryEmitting}
\begin{tabular}{c ccccc ccccc}
\toprule[1pt] \vspace{1pt}
&\multicolumn{3}{c}{40 jobs}&\multicolumn{3}{c}{50 jobs}\\\cmidrule(lr){2-4}\cmidrule(lr){5-7}
Cluster radius: &CR=25&CR=50&CR=100&CR=25&CR=50&CR=100\\\midrule
Mobile coverage   & \multicolumn{1}{c}{6.12} & \multicolumn{1}{c}{6.16}  & 6.18 & \multicolumn{1}{c}{8.04}  & \multicolumn{1}{c}{8.25}  & 8.56  \\ 
\multicolumn{1}{c}{Fixed coverage} & \multicolumn{1}{c}{7.08} & \multicolumn{1}{c}{10.24} & \multicolumn{1}{c}{16.89} & \multicolumn{1}{c}{11.31} & \multicolumn{1}{c}{15.22} & \multicolumn{1}{c}{19.35} \\
\bottomrule[1pt]
\end{tabular}
\end{table}

\section{Conclusion}

Motivated by the advent of 5G technologies, this problem defined a joint routing problem to coordinate the operations of mobile emitters that provide local ICT network coverage and the physical operations of mission vehicles that perform jobs on the ground. This problem has applications across mission-critical environments, including military, search-and-rescue, police and ambulance operations. We developed a double column generation algorithm that add path-based variables for mission vehicles and emitting vehicles iteratively, and enforce a linking constraint between both in a master problem. The algorithm leverages primal and dual information to accelerate the algorithm.

We implemented the model and algorithm in a realistic experimental setup developed in collaboration with Singapore's Defence Science and Technology Agency. Results highlight the benefits of the double column generation algorithm toward solving large-scale instances of the problem, along with the benefits of the optimization methodology developed in this paper to support emerging ICT-physical operations in mission-critical environments. Moreover, the optimization results provide strategic insights regarding the development and deployment of 5G technologies, notably by underscoring the operational benefits of mobile emitters over stationary ones.

This research can be extended in several ways. Notably, an interesting research question involves defining an adaptive mesh discretization as a tractable and high-fidelity approximation of the operations of emitting vehicles in a continuous space. Another extension of the model would involve modeling coverage dissipation as a function of the distance from the 5G emitter, as opposed to assuming a binary coverage variable. Along the same lines, the model of network coverage could be augmented to capture interdependencies between users---including internal users corresponding to multiple mobile emitters in the model, as well as external users in a smart city environment for instance. The models and algorithms developed in this paper provide methodological foundations and practical insights to tackle this emerging class of problems in mission-critical environments.

\bibliographystyle{informs2014}
\bibliography{refs}

\begin{APPENDICES}

\section{Proof of statements}

\subsection{Proof of Proposition~\ref{prop:formulations}}

We draw an equivalence between any mission route $q \in \calQ$ and a collection of $z_{ua}$ variables, as well as any emitting route $p \in \calP$ and a collection of $x_{va}$ variables. Define $\calQ$ as a collection of routes $q$. Each $q$ can be represented as a collection of variables $z_{ua}$ for which $z_{ua} = 1$ for a given $u$, and constraints (2)-(4) (constraints for mission vehicles) and (9)-(13) (time window constraints) are followed. Specifically, these $z_{ua}=1$ variables can be represented as a set of nodes $(i_1, t_1), \ \dots, \ (i_n, t_n)$ in time order. This collection forms a route because of conservation of flow as given by (4). Satisfying the other seven constraints is already embedded within the routes themselves. This definition of $q$ allows for flexibility of time due to the inclusion of time window constraints: when jobs are performed is included in the $y$-variables. To complete the other side of the bijection, any collection of variables $z_{ua}$ satisfying the 8 constraints listed above is a route in $\calQ$.

Similarly, we define $\calP$ as a collection of routes $p$, representable as a collection of variables $x_{va}$ for which $x_{va} = 1$ for a given $v$. In this case, the constraints for emitting vehicles, namely (6), (7), and (8), must apply. As with the mission case, any collection of variables $x_{va}$ satisfying these three constraints is a route in $\calP$.

Now, we show that the constraints satisfied by $q$ and $p$ routes are identical, whether in the explicit or set-partitioning formulations. When the routes are created for set-partitioning, the time windows are already built in to respect all of the time window constraints from (9) to (13); so are constraints (2) to (4), dictating flow constraint and that the route must go out of and back to the depot once. The creation of an emitting route also guarantees that the three constraints from (6) and (8) are satisfied. The linking constraints (14) and (15) are identical to (21), stating that work can only be done if coverage is available. Constraint (5) in the explicit formulation and constraint (18) in the set-partitioning formulation both state that for each job, exactly one vehicle goes into that job. As for set-partitioning constraints (19) and (20), those state that the number of non-trivial routes cannot exceed the number of vehicles allowed. In the explicit formulation, non-trivial routes are still included as routes from the start-depot $(0, 0)$ to the end-depot $(I+1, T+1)$. Because there exists a bijection between any routes in the explicit and set-partitioning formulations, and the objective functions are the same, we conclude that the formulations have the same optimum.

Finally, I claim that the set-partitioning formula is a tighter relaxation than the explicit formulation because we can transform a fractional set-partitioning solution into an explicit formulation, but the other way around is not guaranteed. The explicit formulation contains $y$-variables marking when a job is done. These variables could be fractional. For example, a job with work load 4 and time window size 5 could be done by a mission vehicle with fractional $y$-variables of $0.8$ in the relaxed explicit formulation. However, the set-partitioning formulation does not contain $y$-variables owing to Assumption 1: that jobs are done in the shortest time possible. Thus, a route with fractional $y$-variables permissible under the relaxed explicit formulation is not in the polytope of the relaxed set-partitioning formulation.

\subsection{Proof of Proposition~\ref{prop:workwhenstay}}

In the generation algorithm for feasible mission routes $q \in \calQ$, each route $q$ is defined by an ordered list of nodes $(i, t)$ specifying where the vehicle is at time $t$. Suppose a certain route $q'$ includes some set of nodes at a job location $i'$, $(i', t_1), \ (i', t_2), \ \dots, \ (i', t_{t'})$ but the job is being performed for less than $t'$ time; in other words, there is time slack at the job location. Given the continuity of work assumption, there exists a consecutive block of time $[t_a, t_b]$ over which the job is performed. Then, we can ``redistribute" the time slack before and/or after the job is actually being done into the travel time, so that traveling will have time slack. Whereas the end node of the traveling into job $i'$ was $(i', t_1)$, it will now be $(i', t_a)$, and the start node of the traveling out of job $i'$ will now be $(i', t_b)$ instead of $(i', t_{t'})$. (Note that up to one of these two nodes could be the same if $t_a = t_1$ or $t_b = t_{t'}$.)

\subsection{Proof of Proposition~\ref{prop:earliestWorks}}

Suppose there was a solution $S^{'}$ claimed to be optimal. Consider shifting all jobs to be done at their earliest possible time, and label this solution $S^*$. Then $S^*$ does no worse than $S^{'}$, and perhaps there is an avenue for $S^*$ itself to be revised to do better, specifically, if shifting all the jobs to their earliest possible time enables a certain mission vehicle to take on more jobs that it couldn't have in $S^{'}$, obviating another mission vehicle from having to be deployed. 

In a simplified example, consider two jobs with time windows [4, 10] and [7, 14] and work loads 2 and 5, respectively. Also let travel time between the two be 1. Suppose $S^{'}$ takes one vehicle $u_1$ that goes to the first job from [6, 8] and another vehicle $u_2$ that goes to the second job from [7, 12], then this configuration is suboptimal. Instead, applying early job completion and obtaining solution $S^*$, $u_1$ can work the first job from [4, 6], travel to the second job, and operate it from [7, 12], obviating deployment of $u_2$. This is necessarily better for the objective due to the Triangle Inequality. Letting subscripts be $D$ for depot, $1$ for job 1, and $2$ for job 2, $S^{'}$ has distances $d_{D1} + d_{1D} + d_{D2} + d_{2D}$, while $S^*$ has distances $d_{D1} + d_{12} + d_{2D}$. We see that $d_{1D} + d_{D2} \geq d_{12}$, finishing the proof.

\subsection{Proof of Proposition~\ref{prop:inputBased}}

Proof of Correctness: Suppose an emitting vehicle $v$ arrives or departs before or after one of these constructed times. In the case of arriving before a constructed time $t_c$, then we can "squeeze" the arrival time more tightly into $t_c$ without affecting whether jobs of interest are covered. Similarly, if $v$ departs after a constructed time $t_e$, then we can make it depart earlier at $t_e$ instead without affecting coverage, because constructed times already encapsulate when jobs are being done. By fitting the emitting vehicle's arrival or departure time, we do not affect coverage, but we reduce the number of possible times to consider. Meanwhile, if a emitting vehicle $v$ enters during the middle of any job $i$, there must have been another emitting vehicle $v'$ already present (otherwise, $i$ would not have received coverage, and that solution without $v'$ at the spot would be infeasible). However, having $v$ cover a spot $i$ is already at is redundant and unnecessarily adds to the cost. $v$ would be better deployed at some other uncovered spot. 

\subsection{Proof of Proposition~\ref{prop:dualBased}}

Observe that any considered subproblem emitting route $R$ without any $\xi_{it} > 0$ cannot have a negative reduced cost, because the reduced cost of that route is formed by positive $C^p$, nonnegative $\beta$ (dual variables are nonnegative because the coefficients and decision variables are nonnegative), and a subtraction of a summation relying on multiplication by $\xi_{it}$. If all $\xi_{it} = 0$ for $R$, then the reduced cost $C^p + \beta - \sum \cdot > 0$, which makes it useless. On the contrary, if at least one $\xi_{it} > 0$ is present somewhere along $R$, then there is a possibility that the reduced cost is negative. As no further general filtering can be done on this condition, the solution is to loop over coverage intervals after using input-based pruning which contain a positive $\xi_{it}$. However, not all such intervals are worth traversing, because if there are two intervals $[t_1, t_2]$ and $[t_3, t_4]$ which encompass the same set of $\xi_{it} > 0$, then we can take their intersection to obtain the smallest such interval. It is not hard to see from set theory that we are guaranteed a unique smallest interval, and that it is the ideal interval (because if we take a larger time interval, we risk missing out on going to other intervals for which there exists a $\xi_{it} > 0$).

\subsection{Mission Vehicles Subproblem}\label{appendix:spa}
Please see next page for pseudocode.

\begin{algorithm}[h!]
\floatname{algorithm}{Algorithm}
\caption{Mission Vehicles Subproblem.}
\renewcommand{\arraystretch}{1.0}\small
\label{alg:subproblem}
\begin{algorithmic}[1]
\Require{Number of jobs $I$, travel time array $TT^M$, job distances matrix $d^M$, time windows $T^S$ and $T^E$, work loads $\tau$, dual costs $\rho$, $\pi_i$, $\xi_{it}$.}
\State Initialize paths $\mathbb{S}=\{0\}$, time $T=\{0\}$, reduced cost $RC=\{\rho\}$.

\State Create two pointers for the current state ($p_c$) and total state ($p_t$) to 1.
\While{\text{$p_c \leq p_t$}} \Comment{Still have routes left to check}

\State cur\_path = $\mathbb{S}$[$p_c$], cur\_times = $T$[$p_c$]

\If{last(cur\_path) = $I+1$} \Comment{Check that job nodes still need to be added}
\State{$p_c += 1$} \Comment{The current route is already at the end-depot, go to next route}
\EndIf

\For{i in [1, I+1]} \Comment{Proposed new mission vehicle destination}

\If{$i$ not in cur\_path} \Comment{Cannot revisit already-visited location}

\State proposed\_arrival = last(cur\_times) + $TT^M$[last(cur\_path), $i$]

\If {proposed\_arrival $\leq T^E_i$} \Comment{Feasible arrival time for job $i$}

\State arr = max(proposed\_arrival, $T^S_i$) \Comment{Earliest arrival time}

\State proposed\_cost = $RC[p_c]$ + $d^M$[last(cur\_path), $i$] - $\pi_i$ \Comment{Begin new reduced cost calculation}

\For{t in [arr, arr + $\tau_i$} \Comment{Add xi variables to reduced cost}
\State proposed\_cost += $\xi_{it}$
\EndFor

\If {proposed\_cost $< RC[p_c]$} \Comment{Reduced cost is lower than previous label}

\State Update $\mathbb{S}\gets \text{cur\_path}\cup i$; $T\gets\text{cur\_times}\cup\{\text{arr},\text{arr}+\tau_i$; $RC\gets RC\cup\{\text{proposed\_cost}\}$
\State Increment $p_t$ by 1

\EndIf
\EndIf
\EndIf
\EndFor
\State Increment $p_c$ by 1
\EndWhile
\State Return $\mathbb{S}$, $T$, $RC$
\end{algorithmic}
\end{algorithm}

\end{APPENDICES}

\end{document}